\documentclass[english]{article}
\usepackage[T1]{fontenc}
\usepackage[latin9]{inputenc}
\usepackage[letterpaper]{geometry}
\usepackage{float}
\usepackage{amsmath}
\usepackage{graphicx}
\usepackage{color}
\usepackage{epstopdf}
\usepackage{adjustbox}
\usepackage{diagbox}
\usepackage{multirow}
\usepackage{array}
\usepackage{mathrsfs}
\usepackage{amsmath, amssymb, amsthm}
\usepackage{subfigure}
\usepackage{tikz}

\graphicspath{{./Figures/}}

\usepackage{amsfonts} 
\usepackage{color}
\usepackage{adjustbox}
\usepackage{diagbox}
\definecolor{black}{rgb}{0,0,0}

\definecolor{red}{rgb}{1,0,0}

\definecolor{blue}{rgb}{0,0,1}

\usepackage{multirow}

\makeatother

\usepackage{babel}
\usepackage{authblk}
\title{A local-global multiscale mortar mixed finite element method for multiphase transport in heterogeneous media}

\author[1]{Shubin Fu \thanks{Corresponding Author}}
\author[1]{Eric T. Chung}
\affil[1]{Department of Mathematics, The Chinese University of Hong Kong, Hong Kong SAR}

\begin{document}
		\maketitle
\begin{abstract}
	
	In this paper, we propose a local-global multiscale mortar mixed finite element method (MMMFEM) for multiphase transport in heterogeneous media. 
    We consider the two-phase flow system, the pressure equation is solved
    via the multiscale mortar mixed finite element method, a mass conservation
    velocity field can be obtained, then we use explicit finite volume method
    to solve the saturation equation. We use  polynomials and multiscale basis 
	to form the coarse mortar space. The multiscale basis is the restriction
	of global pressure field obtained at previous time step on the coarse interface. We solve the pressure
	equation on the fine grid to initialize the simulation.
	Numerical experiments on some benchmark 2D and 3D heterogeneous models are
	provided to validate the performance of our method.
\end{abstract}

\section{Introduction}
Modeling flow in porous media has many important practical applications,  such as nuclear waste storage, and oil and natural gas production. However, advanced modern reservoir characterization techniques
can generate very detailed reservoir models with multiple scales, from core
scales (centimeters) to geological scales (kilometers).
Therefore, direct simulation is inevitably expensive and model reduction techniques such as
upscaling and multiscale methods are quite necessary. 
The basis idea of upscaling \cite{durfolsky1991homo,wu2002analysis} is to homogenize the media properties based on some rules so that the problem can be solved in a reduced model.
In multiscale methods \cite{egw10,efendiev2009multiscale,Arbogast_two_scale_04,chung2015mixed,chen2003mixed,jennylt03,Wheeler_mortar_MS_12, ArPeWY07,MHM1,MHM2}, one solves the problems on a coarse grid but with
precomputed multiscale basis functions that includes small scale information of the
media.

Among different multiscale methods, the mixed multiscale finite element method \cite{chen2003mixed}
and the multiscale finite volume method \cite{jennylt03} are especially useful for the 
flow simulation since these methods are mass conservative.
In \cite{aarnes04}, the author used the restriction of
global velocity field in coarse interface as boundary condition
 to compute the multiscale basis functions. This idea
(using limited global information)
is extended to the framework of finite volume \cite{efendiev2006accurate}. 
In \cite{local_global1,local_global2,local_global3},  improvements such as 
adaptive update the multiscale basis functions are developed.
Our goal here is to apply above mentioned ideas to the mortar mixed finite element method (MMFEM)
\cite{mixed_dd1988,arbogast2000}, which is a modification of the mixed finite element methods by introducing a Lagrange multiplier to impose the continuity of flux. 
This modification brings some additional advantages, such as more efficient implementation 
since it yields a a symmetric and positive definite bilinear form that defined only on the
interfaces of the grid. Besides, it also
allows nonconforming grid discretization while keeping the mass conservation property.

The multiscale mortar mixed finite element method (MMMFEM) was proposed in \cite{ArPeWY07,arbogast2013ms} 
to reduce the dimension of the mortar mixed finite element method with the sacrifices of some accuracy.
 In MMMFEM, the construction of the mortar space that defined on coarse interfaces is cubical.
  Polynomials \cite{ArPeWY07} or multiscale  basis functions
  \cite{arbogast2013ms, mortaroffline}
are used to form the mortar space. However, polynomials are only sufficient for very smooth
media,  while multiscale functions lack global information which fails to generate satisfiable
flux for accurate simulation
of coupled flow and transport problems if 
the media is highly heterogeneous\cite{xiao2013multiscale}. 
Developing  efficient  domain decomposition preconditioners \cite{arbogast2013ms} to solve the fine scale problem is a potential way to tackle
this kind of problem. However, if one can equip the mortar space
with basis functions  that can capture global information with the evolving of time, then we can avoid solving  the fine problems.
In \cite{mortaronline}, the authors developed residual based online multiscale basis functions
which contain global information
and successfully applied this method to two phase simulations, however, in some cases, one
needs to iterate several times to get accurate solution, which makes the mortar space 
too enriched. Therefore, further reduction is necessary.
We adopt the above mentioned ideas of using limited global information to 
mortar finite element.
However, naive extension, using the restriction of global pressure on coarse interface to 
form the mortar space fails. Our strategy here is fill the mortar space with fixed 
local basis functions such as polynomials and the restriction of global pressure field obtained at previous time step on the coarse interface.
Therefore, our method is a local-global method. To initialize the simulation, we first need to compute an accurate 
initial global pressure field, which can be achieved by solving the fine problem with
preconditioners \cite{arbogast2013ms} or online multiscale method \cite{mortaronline}.
Numerical results show that with only 2 bases in 2D case and 4 bases in 3D cases,
the MMMFEM can generate very accurate saturation profiles with much reduced computational 
cost.

The paper is organized as follows. In Section \ref{sec:two-phase}, we introduce
the coupled two-phase flow system.  In Section \ref{sec:multiscale}, we present
 MMMFEM method for the pressure equation. The construction of the local-global multiscale mortar space in
introduced in Section \ref{sec:basis}.  Numerical examples are given in Section \ref{numerical-results}, and conclusions are made in the last section.

\section{Mathematical model}\label{sec:two-phase}

We consider  a model problem for immiscible and
incompressible two-phase flow in a reservoir domain (denoted by $\Omega$)
with the assumption that the fluid displacement is driven by viscous effects, that is, we neglect compressibility
and gravity for simplicity. 
We consider water and oil phases which are assumed to be immiscible.
By the mass conservation law, we have following equations for each phase
$l$:
\begin{equation} \label{eq:mass}
\phi\frac{\partial S_l}{\partial t} + \nabla \cdot {\bf u}_l = q_l,
\end{equation}
where the phase velocities ${\bf u}_l$ are given by Darcy's law

\begin{equation} \label{eq:darcy}
{ \bf u}_l=-\lambda_l(S_l) { K} \nabla {p_l}.
\end{equation}
Here  $\phi$ is the porosity, $S_l$ is the $l$-phase saturation (fraction of the void occupied by
phase $l$), and $q_l$ is a source (or sink) term,
 ${K}$ is the permeability
 tensor which can be highly heterogeneous, $p_l$ is the phase pressure, and $\lambda_l$ is the phase mobility given by 
 $\lambda_l(S_l)=\frac{k_{rl}(S_l)}{\mu_l}$ , where $k_{rl}$ and $\mu_l$ are the relative permeability and viscosity of phase
 $l$, respectively. 

 Since we neglect effects from
 capillary pressure so that $\nabla p_o = \nabla p_w$, and we denote $p_o = p_w = p$. Then we can
 derive the following {\it pressure equation}
\begin{eqnarray}
\bf u &=& -(\lambda_w + \lambda_o) \nabla p, \\
\nabla\cdot \bf u &=& q,		
\end{eqnarray}
where $\bf u=\bf u_o+\bf u_w$ and $q=q_o+q_w$

We assume that the two phases occupy
the void space completely, therefore $S_o+S_w=1$, let $\lambda(S)=\lambda_w(S_w)+\lambda_o(1-S_w ) $ be the total mobility and the  water {fractional flow} $f_w = \lambda_w /\lambda$,  then we have the conservation
equation for water as follows
\begin{equation} \label{eq:sat}
\phi\frac{\partial S_w}{\partial t} + \nabla \cdot ({f_w\bf u}_l) = q_w,
\end{equation}
which is henceforth called the {\it saturation equation},
In the following we will, drop the w-subscripts of $S_w$ for ease of notation, 
We assume no flow at boundary and the initial saturation is given.
Note the pressure equation and the saturation equation   are coupled through the total mobility.
Here, we follow the sequential formulation \cite{aarnes04} to solve the above coupled system, that is we use saturation at the previous time
step to solve the pressure equation to obtain the velocity $\bf u$. 
With the updated velocity, we solve the transport equation with explicit finite volume method and obtain
new saturation, which is then used to solve the pressure equation again and the process is continued
until the final time is reached.
\section{Multiscale mortar mixed finite element method}\label{sec:multiscale}

In this section, we will introduce the framework of multiscale mortar mixed finite element method
for the pressure equation.
To fix our attention,
We consider
\begin{subequations}\label{eq:problem}
	\begin{alignat}{2}
	\label{original equation-1}
	\boldsymbol{u} + \lambda {K}\nabla p &= 0 \qquad && \text{in $\Omega$,}\\
	\label{original equation-2}
	\nabla \cdot \boldsymbol{u}  &= q && \text{in $\Omega$,}\\
	\label{boundary condition}
	\boldsymbol{u\cdot n} &= 0 && \text{on $\partial \Omega$,}
	\end{alignat}
\end{subequations}
where $\Omega \subset \mathbb{R}^d (d=2, 3)$ is a bounded polyhedral domain with outward unit normal vector $\boldsymbol{n}$ on the boundary,
$q \in L^2(\Omega)$.
$K(x)$ is a symmetric and uniformly
positive definite permeability tensor with components in $L^{\infty}(\Omega)$.
\subsection{Formulation of the method}
 We divide $\Omega$  into non-overlapping polygonal coarse elements $K_i$ 
so that $\overline{\Omega}=\cup_{i=1}^N\overline{K}_i$,
where $N$ is the number of coarse elements.
The decomposition of the domain can be nonconforming. We call $E_H$ a coarse
interface of the coarse block $K_i$ if $E_H = \partial K_i \cap \partial K_j $ or $E_H= \partial K_i \cap \partial{\Omega}$.
Let $\mathcal{E}_H(K_i)$ be the set of all coarse interfaces on the boundary of the coarse element $K_i$, and $\mathcal{E}_H=\cup_{i=1}^N\mathcal{E}_H(K_i)$
be the set of all coarse interfaces. For a coarse element $K_i$
and coarse interface $E_i$, let
\begin{equation*}
{\bf V}_i=H(\text{div};K_i), \quad{\bf V}=\bigoplus_{i=1}^N{\bf V}_i,
\quad W_i=L_2(K_i),\quad W=\bigoplus_{i=1}^N W_i,
\quad M_i=H^{1/2}(E_i),\quad M=\bigoplus_{i=1}^{|\mathcal{E}_H|} M_i,
\end{equation*}
Then, a weak form of system (\ref{original equation-1})-(\ref{boundary condition}) reads: find
${\bf u}\in {\bf V}$, $p\in W$ and $\Lambda \in M$ such that for each $1\leq i \leq N$,
\begin{subequations}\label{Eq: variation form}
	\begin{alignat}{2}
	\label{variation-1}
	((\lambda {K})^{-1} \boldsymbol{u}, \boldsymbol{v})_{K_i} - (p, \nabla\cdot\boldsymbol{v})_{K_i}+\left\langle\Lambda, \boldsymbol{v}\cdot \boldsymbol{n}_i\right\rangle_{ {\partial K_i}} &=0\quad && \forall ~ \boldsymbol{v}\in \boldsymbol{V}_i,\\
	\label{variation-2}
	(\nabla\cdot \boldsymbol{u}, w)_{K_i}&= (q,w)_{K_i}\quad&& \forall ~ w\in  W_i,\\
	\label{variation-3}
	\sum_{i=1}^N\left\langle \boldsymbol{u}\cdot \boldsymbol{n}_i, \mu\right\rangle_{{\partial K_i}} &= 0\quad && \forall ~ \mu\in M.
	\end{alignat}
\end{subequations}
where $\boldsymbol{n}_i$ is the outer unit normal to $\partial K_i$.
$\Lambda$ is a Lagrange multiplier introduced
to impose the  continuity of $\boldsymbol{v}\cdot \boldsymbol{n}$ on coarse interface.

\subsection{The finite element approximation}
We further partition each coarse element $K_i$ into a finer mesh with mesh size $h_i$.
Let $\mathcal{T}_h=\cup_{i=1}^N\mathcal{T}_h(K_i)$ be the union of all these partitions, which forms a fine partition of the domain $\Omega$.
We use $h=\text{max}_{1\leq i\leq n}h_i$ to denote the mesh size of $\mathcal{T}_h$.
In addition, we let
$\mathcal{E}_h(K_i)$
be the set of all edges (or faces in three dimensions) of the partition $\mathcal{T}_h(K_i)$ 
and let $\mathcal{E}_h=\cup_{i=1}^N\mathcal{E}_h(K_i)$ be the set of all edges (or faces) in the partition $\mathcal{T}_h$. 
Let ${\bf V}_{h,i}\times W_{h,i}\subset {\bf V}_{i}\times W_{i}$ be any of the mixed finite element spaces satisfying the inf-sup condition such as the Raviart-Thomas spaces (\cite{raviart1977mixed}). Define
${\bf V}_h=\oplus_{i=1}^N {\bf V}_{h,i}$ and $W_h=\oplus_{i=1}^N
W_{h,i}/\mathbb{R}$ for the global discrete flux and pressure. Let $M_{H,i}\subset L_2(E_{i})  
$  be the local coarse mortar finite space, and
	$M_H=\oplus_{1\leq i \leq {|\mathcal{E}_H|}}M_{H,i}$  be the entire coarse  mortar finite element space.
Then,  the multiscale mortar mixed finite element approximation of (\ref{variation-1})-(\ref{variation-3}) is to seek
${\bf u}_h\in {\bf V_h}$, $p_h\in W_h$ and $\Lambda_H \in M_H$ such that for each $1\leq i \leq N$,
\begin{subequations}\label{Eq: discrete form}
	\begin{alignat}{2}
	\label{discrete-1}
	((\lambda {K})^{-1} \boldsymbol{u}_h, \boldsymbol{v}_h)_{K_i}-(p_h, \nabla\cdot\boldsymbol{v}_h)_{K_i} +\left\langle\Lambda_H, \boldsymbol{v}_h\cdot \boldsymbol{n}_i\right\rangle_{ {\partial K_i}}&=0 \quad && \forall ~\boldsymbol{v}_h\in \boldsymbol{V}_{h,i},\\
	\label{discrete-2}
	(\nabla\cdot \boldsymbol{u}_h, w_h)_{K_i}&= (q,w_h)_{K_i}\quad&& \forall ~ w_h\in  W_{h,i},\\
	\label{discrete-3}
	\sum_{i=1}^N\left\langle \boldsymbol{u}_h\cdot \boldsymbol{n}_i, \mu_H\right\rangle_{{\partial K_i}} &= 0\quad && \forall ~\mu_H\in M_H.
	\end{alignat}
\end{subequations}
Local conservation within each coarse element can be obtained by setting $w_h=1$ in Equation (\ref{discrete-2}), weak continuity of the flux is enforced across the coarse interfaces via  Equation (\ref{discrete-3}).
To make the above approximation be well posed, the two scales must be chosen such that
 the mortar space is not too rich compared to the normal traces of the 
 velocity spaces in coarse element.
 For the existence of the solution to above system, we refer \cite{ArPeWY07}
and the reference therein.
\subsection{Interface problem}

The main feature of the MMMFEM is that it could be implemented by just solving a small dimension system 
defined on the coarse interface together with the solutions of some local problems in each coarse element.

Define bilinear forms $a_{H,i}: M_H\times M_H \rightarrow \mathbb {R}, i=1,\cdots, N$ by
$$a_{H,i}(\Lambda, \mu)= - \left\langle {\bf u}_h^{\ast}(\Lambda)\cdot {\bf n}_i, \mu \right\rangle_{{\partial K_i}},$$ and $a_H: M_H\times M_H \rightarrow \mathbb {R}$ by
$$a_H=\sum_{i=1}^N a_{H,i}(\Lambda, \mu), $$
where $\big({\bf u}_h^{\ast}(\Lambda), p_h^{\ast}(\Lambda)\big)\in  {\bf V}_{h} \times W_{h}$ solves ($\Lambda$ given, $q=0$)
\begin{subequations}\label{Eq: domain dec1}
	\begin{alignat}{2}
	\label{domain dec1-1}
	\big((\lambda {K})^{-1} \boldsymbol{u}_h^{\ast}(\Lambda), \boldsymbol{v}_h\big)_{K_i}-\big(p_h^{\ast}(\Lambda), \nabla\cdot\boldsymbol{v}_h\big)_{K_i}&= -\left\langle\Lambda, \boldsymbol{v}_h\cdot \boldsymbol{n}_i\right\rangle_{ {\partial K_i}} \qquad && \forall ~ \boldsymbol{v}_h\in \boldsymbol{V}_{h,i},\\
	\label{domain dec1-2}
	\big(\nabla\cdot \boldsymbol{u}_h^{\ast}(\Lambda), w_h\big)_{K_i}&= 0\qquad&& \forall ~ w_h\in  W_{h,i},
	\end{alignat}
\end{subequations}
for each $1\leq i\leq N.$

Define linear functionals $g_{H,i}: M_H \to \mathbb {R}$ by
$$g_{H,i}(\mu)= \left\langle  \boldsymbol{\bar{u}}_h \cdot {\bf n}_i, \mu \right\rangle_{{\partial K_i}},$$ and $g_H: M_H \to \mathbb {R}$ by
$$g_H(\mu)=\sum_{i=1}^N g_{H,i}(\mu),$$
where $(\boldsymbol{\bar{u}}_h, \bar{p}_h)\in  {\bf V}_{h} \times W_{h}$ solves ($\Lambda=0, q$ given) for $1\leq i\leq N$
\begin{subequations}\label{Eq: domain dec2}
	\begin{alignat}{2}
	\label{domain dec2-1}
	\big((\lambda K)^{-1} \boldsymbol{\bar{u}}_h, \boldsymbol{v}_h\big)_{K_i}-\big(\bar{p}_h, \nabla\cdot\boldsymbol{v}_h\big)_{K_i}&= 0 \quad && \forall ~ \boldsymbol{v}_h\in \boldsymbol{V}_{h,i},\\
	\label{domain dec2-2}
	\big(\nabla\cdot \boldsymbol{\bar{u}}_h, w_h\big)_{K_i}&= (q,w_h)_{K_i}\quad&& \forall ~ w_h\in  W_{h,i}.
	\end{alignat}
\end{subequations}

Then, the system \ref{Eq: discrete form} reduce to find $\Lambda_H\in M_H$ such that
\begin{equation}\label{interface}
a_H(\Lambda_H, \mu)=g_H(\mu), \quad \forall ~\mu \in M_H.
\end{equation}
It can be shown (see \cite{arbogast2000}) that the interface problem (\ref{interface})  produces the solution of
(\ref{discrete-1})-(\ref{discrete-3}) via
$$\boldsymbol{u}_h= {\bf u}_h^{\ast}({\Lambda_H})+\boldsymbol{\bar{u}}_h, p_h=\tilde{p}_h-\frac 1 {|\Omega|}\int_{\Omega}\tilde{p}_h, $$
where $\tilde{p}_h=p_h^{\ast}({\Lambda_H})+\bar{p}_h. $

The bilinear form $a_H(\cdot,\cdot)$ defined on interface is symmetric and positive semi-definite on $M_H$.
There are several ways to implement the MMMFEM such as 
using a multiscale flux basis \cite{ganis2009implementation}.
In their implementation, the matrix corresponding to
the bilinear form will not be formed explicitly since they solve the linear system iteratively, after the multiscale flux basis is computed, only 
linear combinations are needed. Here, we use another method, that is, we
assemble the matrix explicitly and apply direct solver to solve the linear
system since its dimension is very small. Parallel computed can be easily adopted for the 
local problems.
We remark that the MMMFEM shares some similarities with the Multiscale Hybrid-Mixed finite element methods 
\cite{MHM1,MHM2}, both methods involve solving some local problems and a small dimension global problem.
However, the local problems in the Multiscale Hybrid-Mixed finite element methods are Neumann boundary condition
problems.

\section{Basis space}\label{sec:basis}
The coarse space $M_H$ is quite important
to the accuracy and efficiency of the MMMFEM.
If the mortar space is small, then the resulting linear
system can be solved  relatively easily, however the accuracy
can not be guaranteed since the continuity of
the flux only be imposed very weakly. If the mortar space is 
too enriched, then the method will still be expensive.
Various types of mortar were proposed in literature.
In \cite{ArPeWY07}, the authors considered the polynomials for the 
mortar space and proved some convergence results.
Several media dependent multiscale mortar space are also studied,
for example, in \cite{mortaroffline}, an enriched multiscale space was designed
for the high-contrast flow. Homogenized multiscale basis (see \cite {xiao2013multiscale}) was constructed for the flow in highly heterogeneous
media, which can yields a moderate accurate velocity field. However this 
velocity field is not accurate enough to be used to transport some substance
(see \cite {arbogast2015two}).
Residual driven multiscale online basis 
that includes global information was introduced in \cite{mortaronline},
successful applications for the two-phase flow simulations were presented.
However, for some tough cases, the mortar space may be large ($10\sim15$ local bases are needed in three dimensions). 
Our goal here is to further reduce the size of $M_H$
while keeping accuracy. We adopt the idea of using limited global information
(see \cite{local_global1,local_global2,local_global3,efendiev2006accurate,aarnes})
to construct the mortar space.
However, only use the restriction of the global pressure field
on the coarse interface to form $M_H$ is not enough. In particular, 
for each coarse interface, 
we will use two types functions to fill local $M_{H,i}$ at time step
$t$, polynomial $P_r$ ($r$ is the polynomial order)
and the restriction of $\Lambda_H$ obtained with the MMMFEM 
at time step $t-1$. Therefore, the space is different at different time steps.
 We use the MMFEM to solve the initial pressure equation.
In practice, $r=0$ or $1$ is enough for simulation, therefore the local mortar space
$M_{H,i}$ only consists $2\sim 4$ bases.

\section{Numerical examples}\label{numerical-results}
\begin{figure}[H]
	\centering
	\subfigure[$K_{1}$ in $\log_{10}$ scale, the precision is $60 \times 220$]{
		\includegraphics[trim={0 1cm 0  0cm},clip,width=4in]{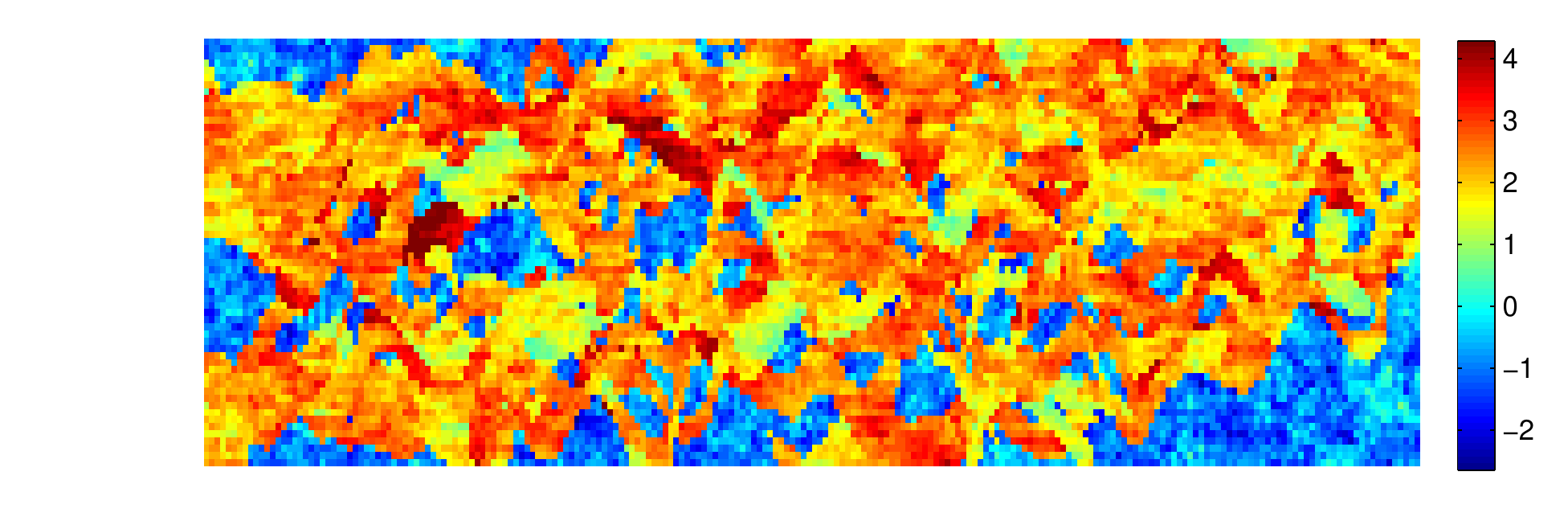}}	
	\subfigure[$K_{2}$ in $\log_{10}$ scale, the precision is $ 60\times220\times30$]{
		\includegraphics[trim={0 1cm 0  1cm},clip,width=3in]{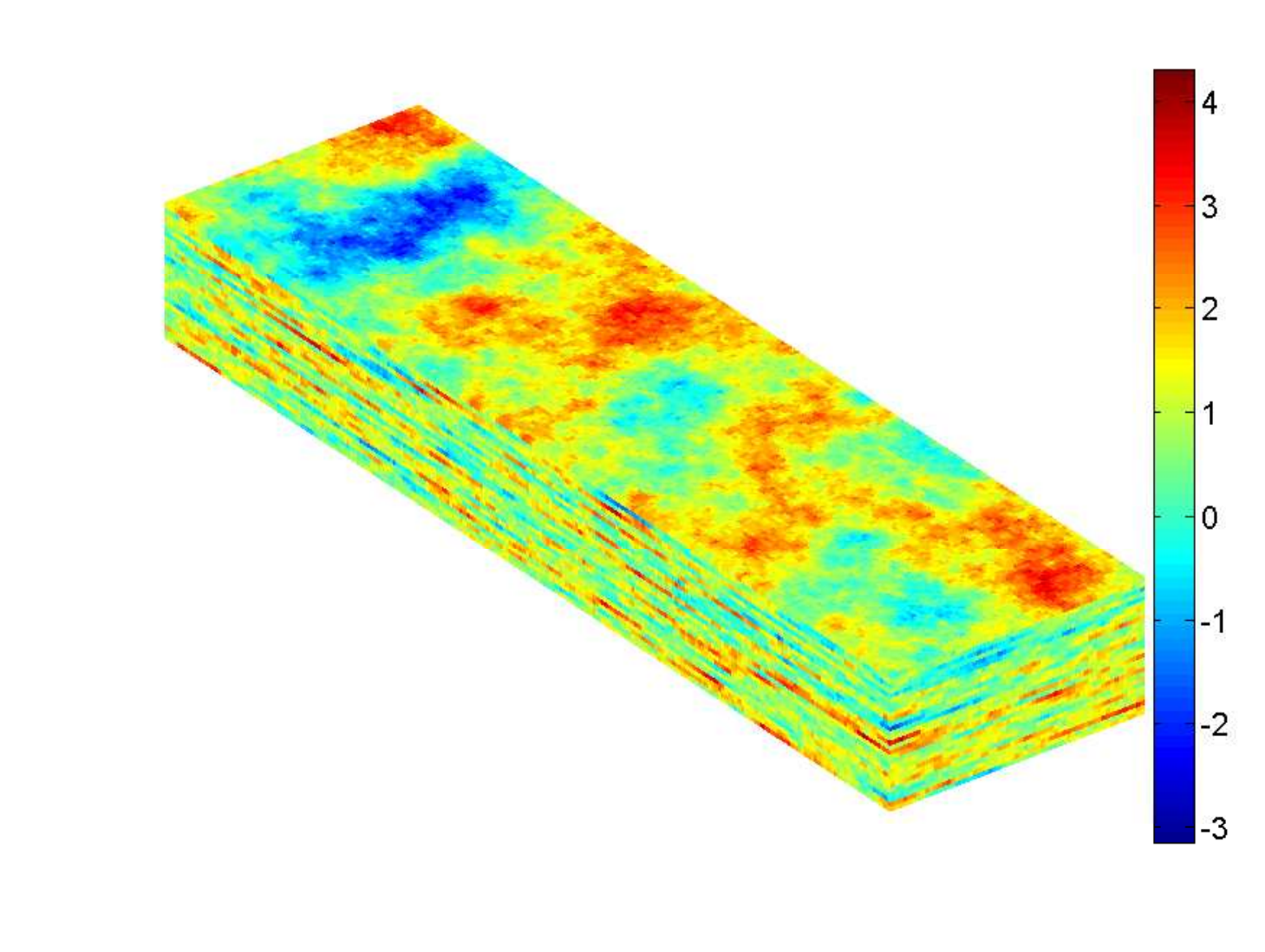}}
	\subfigure[$K_{3}$ in $\log_{10}$ scale, the precision is $ 60\times220\times50$]{
		\includegraphics[trim={0 1cm 0  1cm},clip,width=3in]{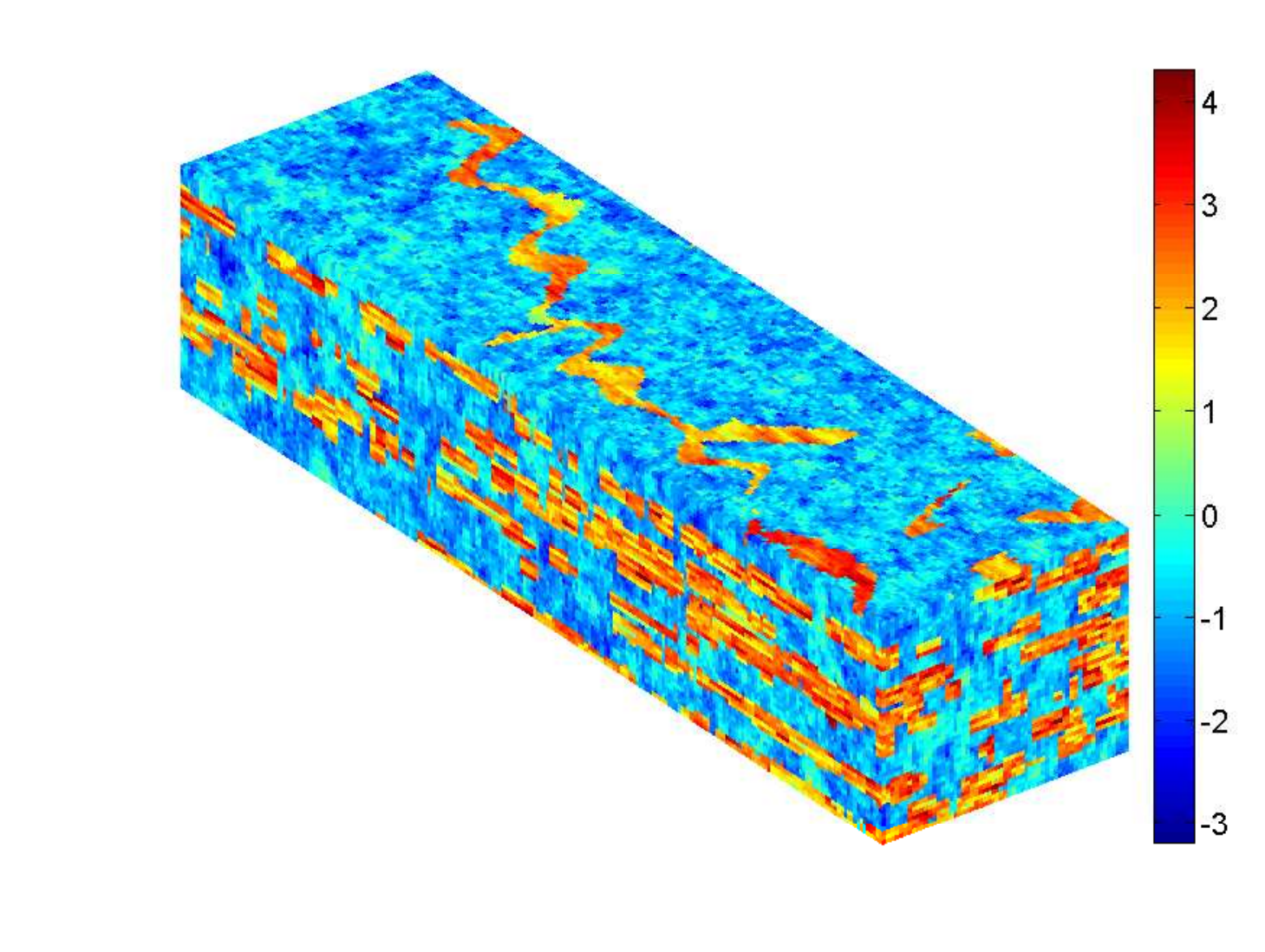}}	
	\caption{Permeability fields.}
	\label{fig:model} 
\end{figure}

In this section, we present several representative numerical experiments to show the performance of our method.
We consider three highly heterogeneous permeability fields $K$ depicted in Figure \ref{fig:model}. 
All these permeability fields are extracted from the tenth SPE comparative solution project (SPE10)  \cite{Aarnes2005257},
which is commonly used benchmark permeability field  to test different upscaling techniques and multiscale methods.
The first model is the last layer of the SPE10 dataset, as it is shown, it is highly heterogeneous and contains long channels.
 The second model is the first 30  layers of the SPE10 dataset,
 which represents a prograding near-shore environment.
The third model is extracted from the last 50 layers of the SPE10 dataset
which represents an Upper Ness fluvial formation, with clearly 
visible channels. The model 2 is smoother than the model 3. The precision of these models are $60\times220$,
$60\times220\times30$ and $60\times220\times50$ respectively.
We will consider two coarsening ratio, i.e, $n=5$ and $n=10$. 
In the reported tables below, 
we use "$Nb$" to represent the number of basis functions used 
in the local mortar space, and "Dim" to represent the dimension
of the resulting coarse linear system of the MMMFEM.

We use the quadratic relative permeability curves $k_{ro} = ( 1-S ) ^2$ and $k_{rw} = S^2$ and, along with $\mu_w = 1$ and
$\mu_o = 5$ for the fluid viscosities. 
We set the initial water saturation as zero in the whole reservoir.  
We include 4 injectors in the corners of the reservoir and 1 producer
in the center of the reservoir, we will record the oil production in this producer. 
The water cut (or fractional flow) is defined as the fraction of water in the produced fluid and is
given by $\frac{q_w}{q_t}$, where $q_t = q_o + q_w$, with $q_o$ and $q_w$ being the flow rates of oil and water at the producer.
We apply 
10 times cheap Jacobi  iterations \cite{mansfield1991damped} to smooth the multiscale  solution.
The total simulation time is 2000, we update the 
velocity at every  50 time steps, i.e,. the pressure equation is solved
40 times. 
 We define the relative saturation at time step $i$ $(i=1,2\cdots,40)$ as\begin{equation*}
e_S(i):=\frac{\|S_{{ms}}(i)-S_{ref}(i)\|_{(L^2,\Omega)}}{\|S_{ref}(i)\|_{(L^2,\Omega)}}.
\end{equation*}
We also denote $e_S$ as the average saturation error history over the time instants.
For the 2D problem, we will consider two mortar spaces, i.e, multiscale basis and multiscale 
basis together with constant function (P0). For the 3D problems, we will consider using multiscale
basis with constant function or first order polynomials (P1).
Since  two coarsening ratio cases will be tested, therefore, 4 multiscale mortar solutions will be computed for all test models. The reference solution is computed with the MMFEM on the fine grid.

\subsection{2D test}
\begin{figure}[H]
	\centering
		\includegraphics[trim={0 1.1cm 0  0cm},clip,width=4in]{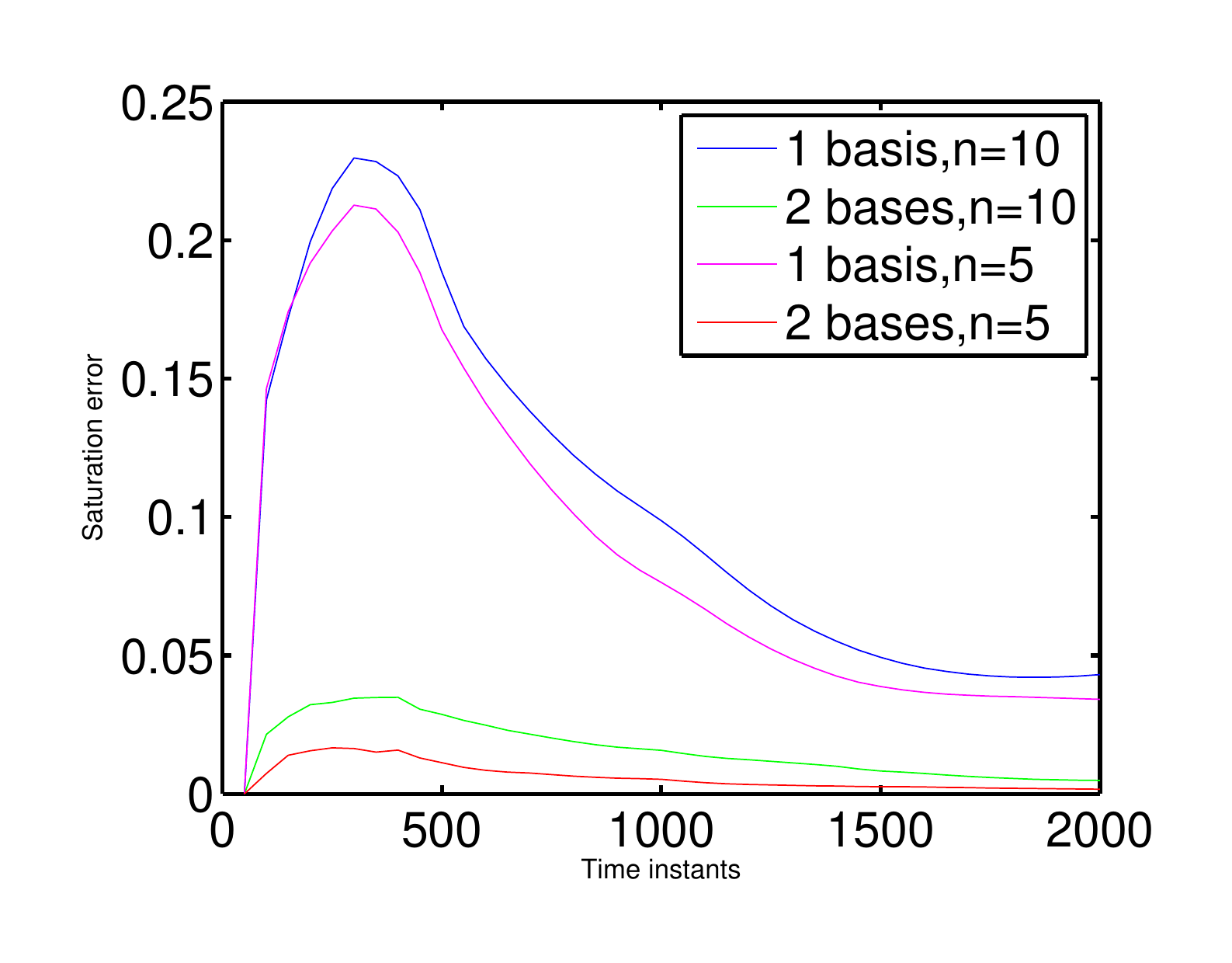}	
	\caption{Saturation error history for model 1.}
	\label{fig:spe85saterr} 
\end{figure}

\begin{figure}[H]
	\centering
		\includegraphics[trim={0 1.1cm 0  0cm},clip,width=4in]{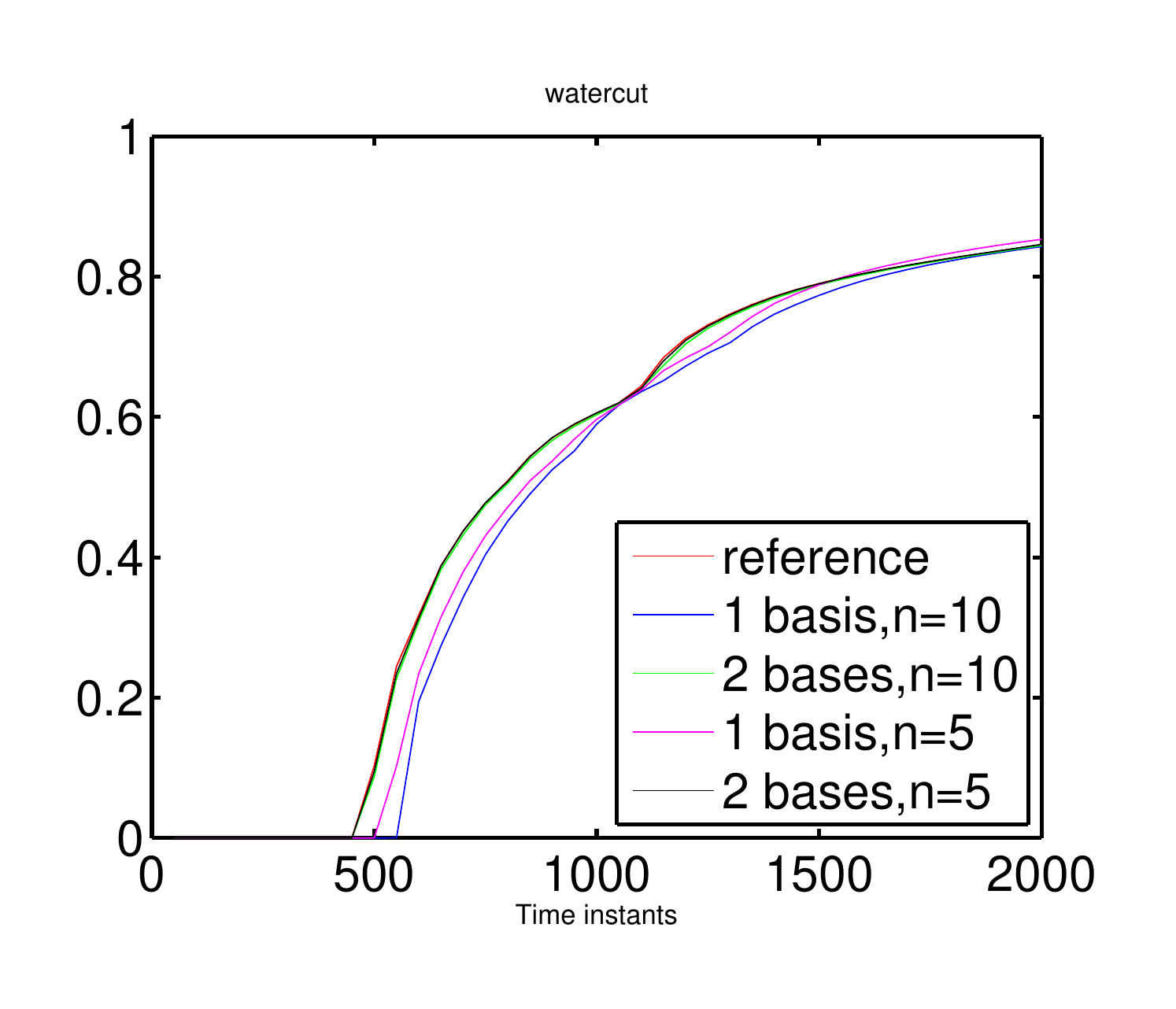}		
	\caption{Watercut comparison for model 1.}
	\label{fig:spe85pc} 
\end{figure}

\begin{figure}[H]
	\centering
	\subfigure[Reference saturation profile at $t = 1000$]{
		\includegraphics[trim={0 1cm 0  0.8cm},clip,width=4in]{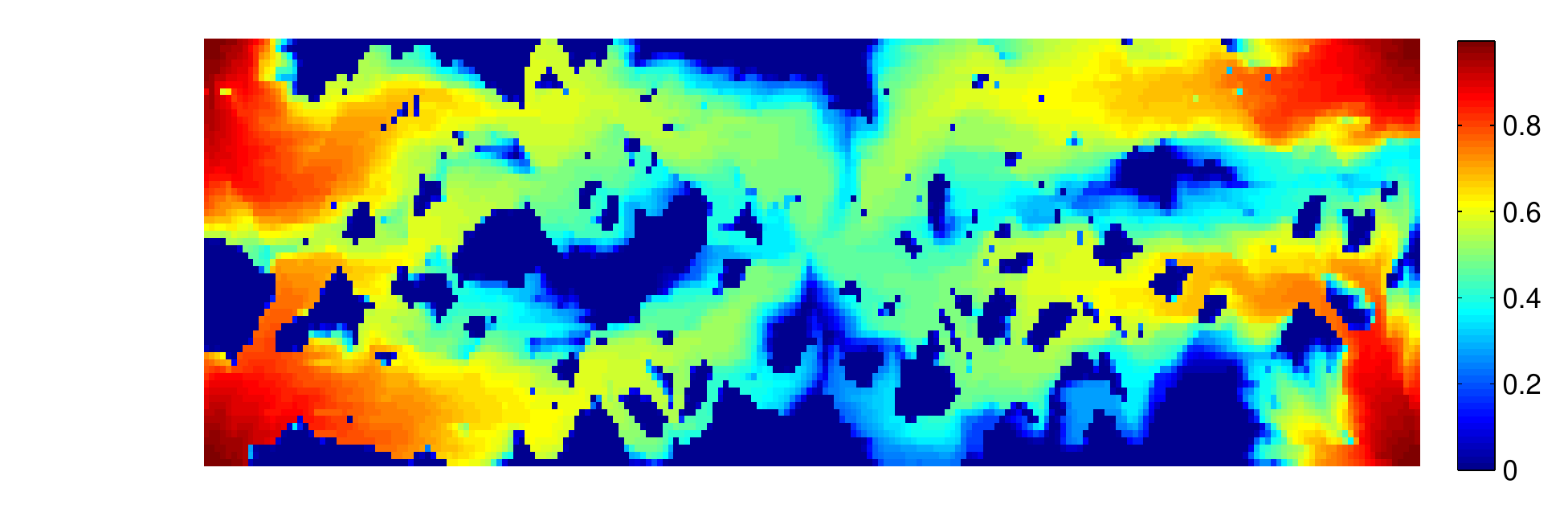}}
	\subfigure[Saturation profile with MMMFEM at $t = 1000$, $n=10$, $Nb=1$]{
		\includegraphics[trim={0 1cm 0  0.8cm},clip,width=4in]{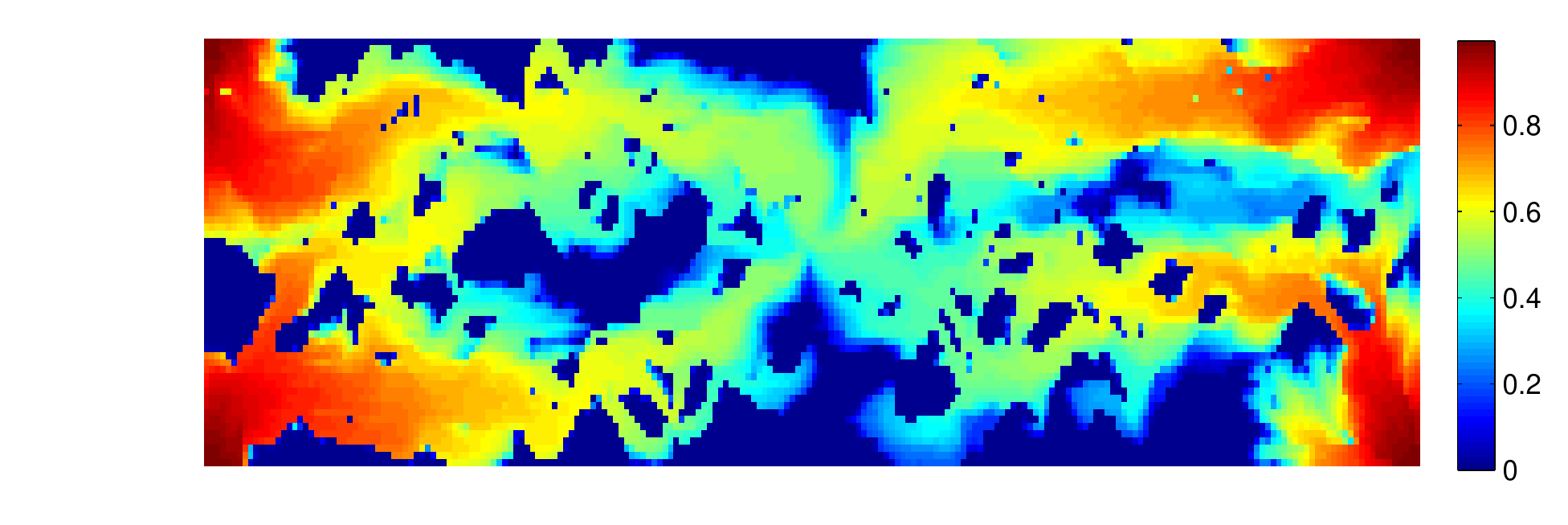}}
	\subfigure[Saturation profile with MMMFEM at $t = 1000$, $n=5$, $Nb=2$]{
		\includegraphics[trim={0 1cm 0  0.8cm},clip,width=4in]{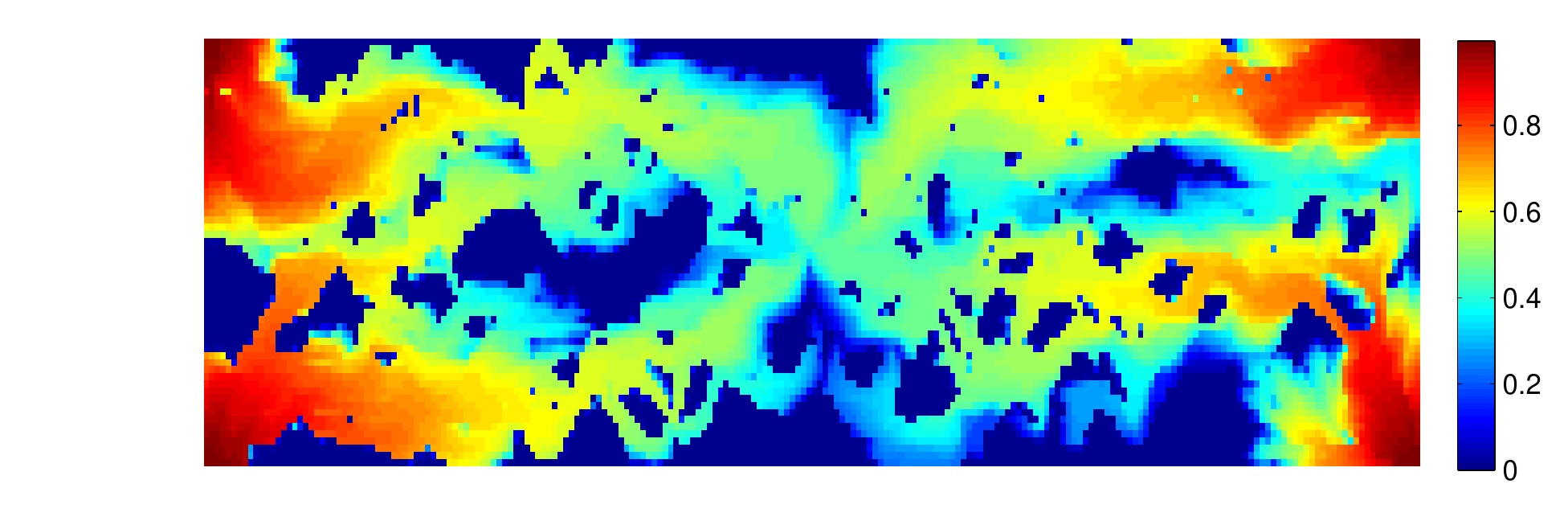}}	
	\caption{Saturation profile at $t = 1000$ for fine-scale solution (top figure), multiscale mortar solution when $n=10$ and $Nb=1$ (middle figure), and multiscale mortar solution when $n=5$ and $Nb=2$ (bottom figure).}
	\label{fig:sat85} 
\end{figure}

Figure \ref{fig:spe85saterr} shows the saturation error 
history of the multiscale
solution with different test setting with respect to the
reference solution  for test model 1. We can see that 
if only the multiscale basis is used, the error is large 
especially at earlier times.
if we add constant functions to the mortar space, the
accuracy of the multiscale solution is improved significantly,
In Figure \ref{fig:spe85pc}, we plot the watercut comparison
between different multiscale solution and the reference solution
for test model 1. We also observe that if two basis functions 
are used, the multiscale solution can approximate the reference
solution pretty well. Actually, it is hard to distinguish 
the line that represents the multiscale watercut when $Nb=2$ and $n=5$  and the reference watercut line.
Figure \ref{fig:sat85} shows the saturation profile comparison between
the multiscale solution and the reference solution.
We find that when only 1 basis is used and $n=10$, the saturation profile 
(middle figure in \ref{fig:sat85})
generated with the MMMFEM  fails to capture all the details of 
reference saturation profile (top figure in Figure \ref{fig:sat85}), however if we use 2 bases and increase
the coarse-grid size, the MMMFEM solution (bottom figure in Figure \ref{fig:sat85}) is almost identical as the reference solution.

\subsection{3D tests}

\begin{figure}[H]
	\centering
	\centering
	\includegraphics[trim={0 .2cm 0  0cm},clip,width=4in]{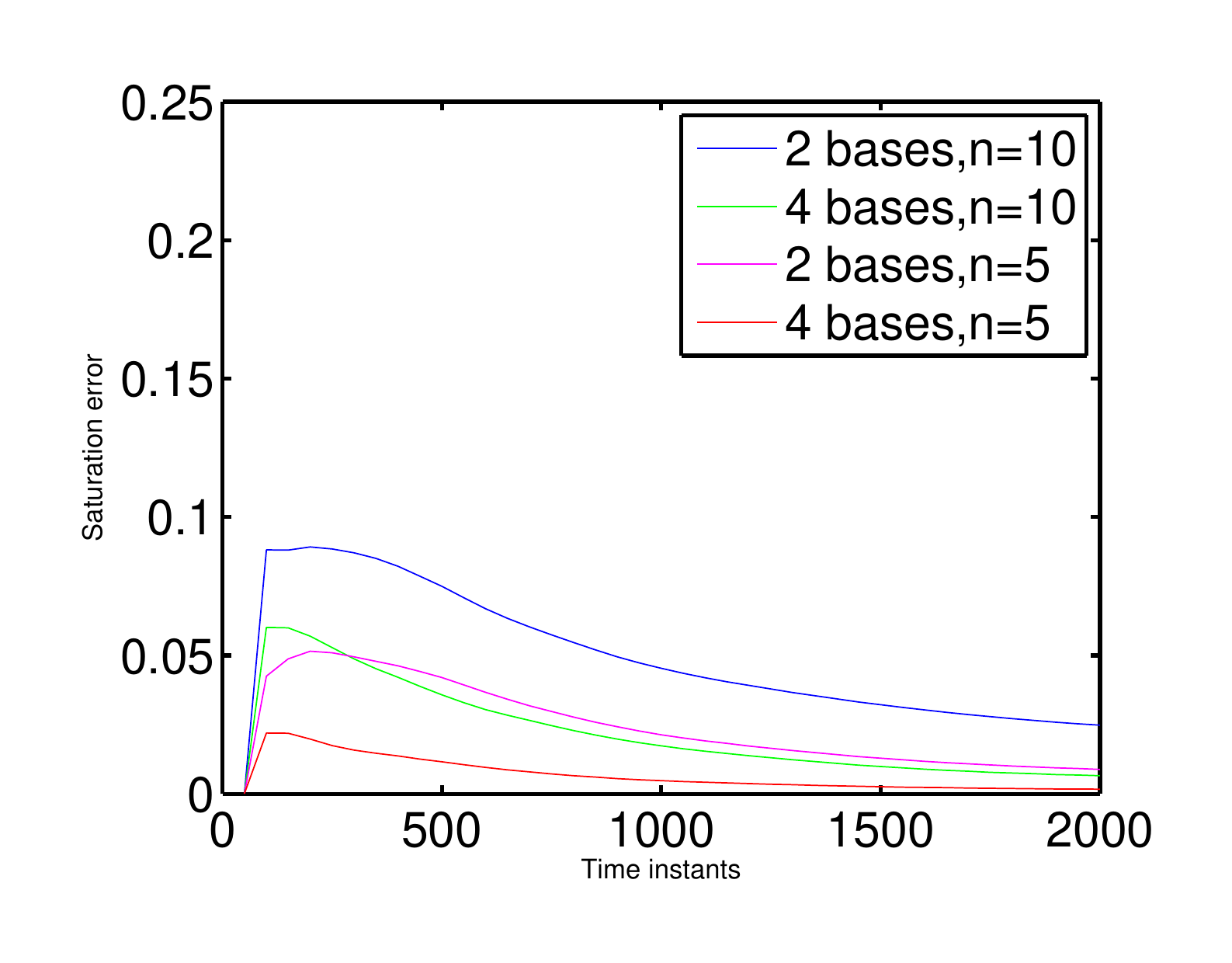}	
	\caption{Saturation error history for model 2.}
	\label{fig:satfirst30} 
\end{figure}

\begin{figure}[H]
	\centering
	\includegraphics[trim={0 .2cm 0  0cm},clip,width=4in]{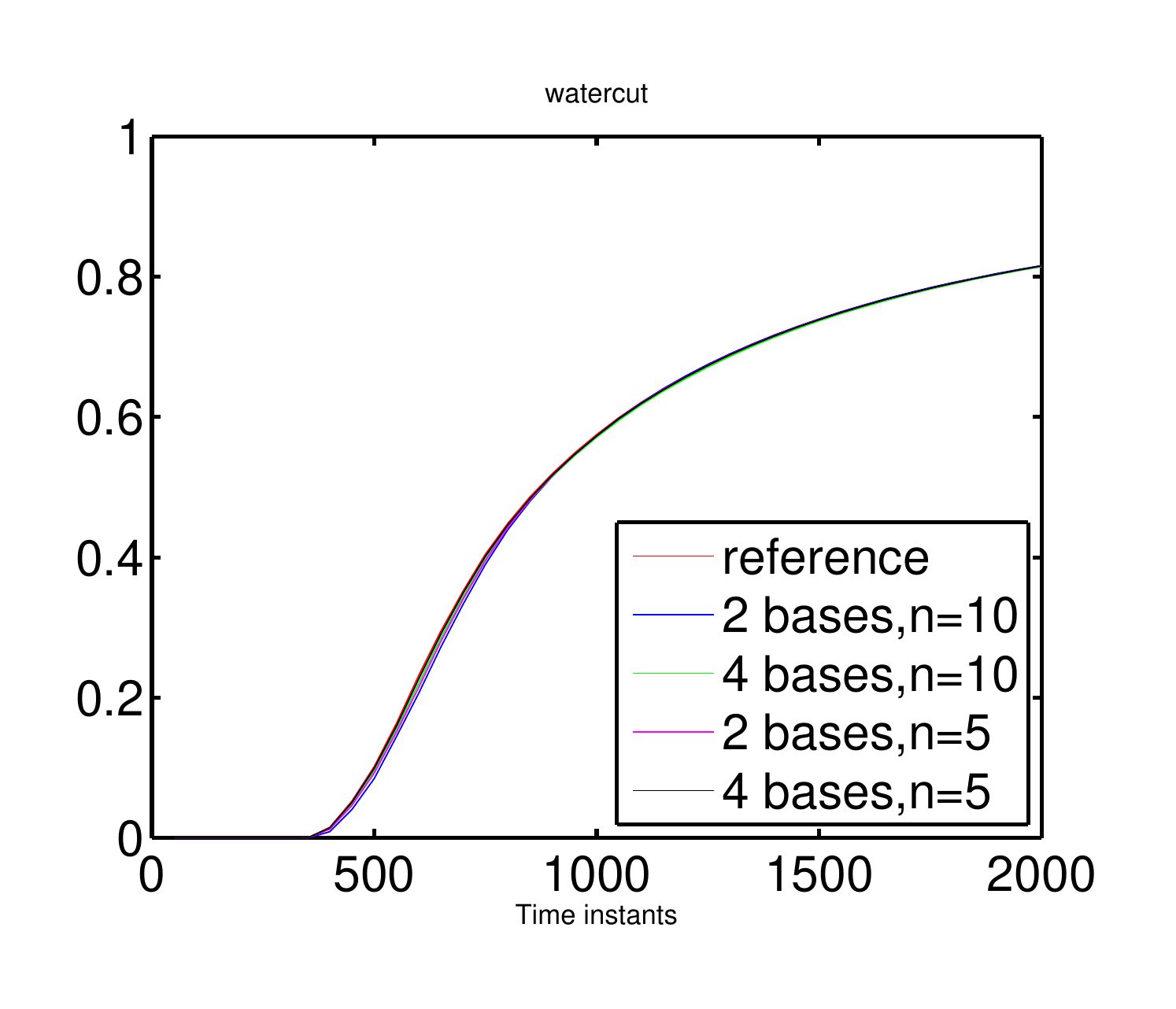}		
	\caption{Watercut comparison for model 2.}
	\label{fig:pcfirst30} 
\end{figure}

\begin{figure}[H]
	\centering
	\subfigure[Reference saturation profile at $t = 1000$]{
		\includegraphics[trim={0 1cm 0  0.8cm},clip,width=4in]{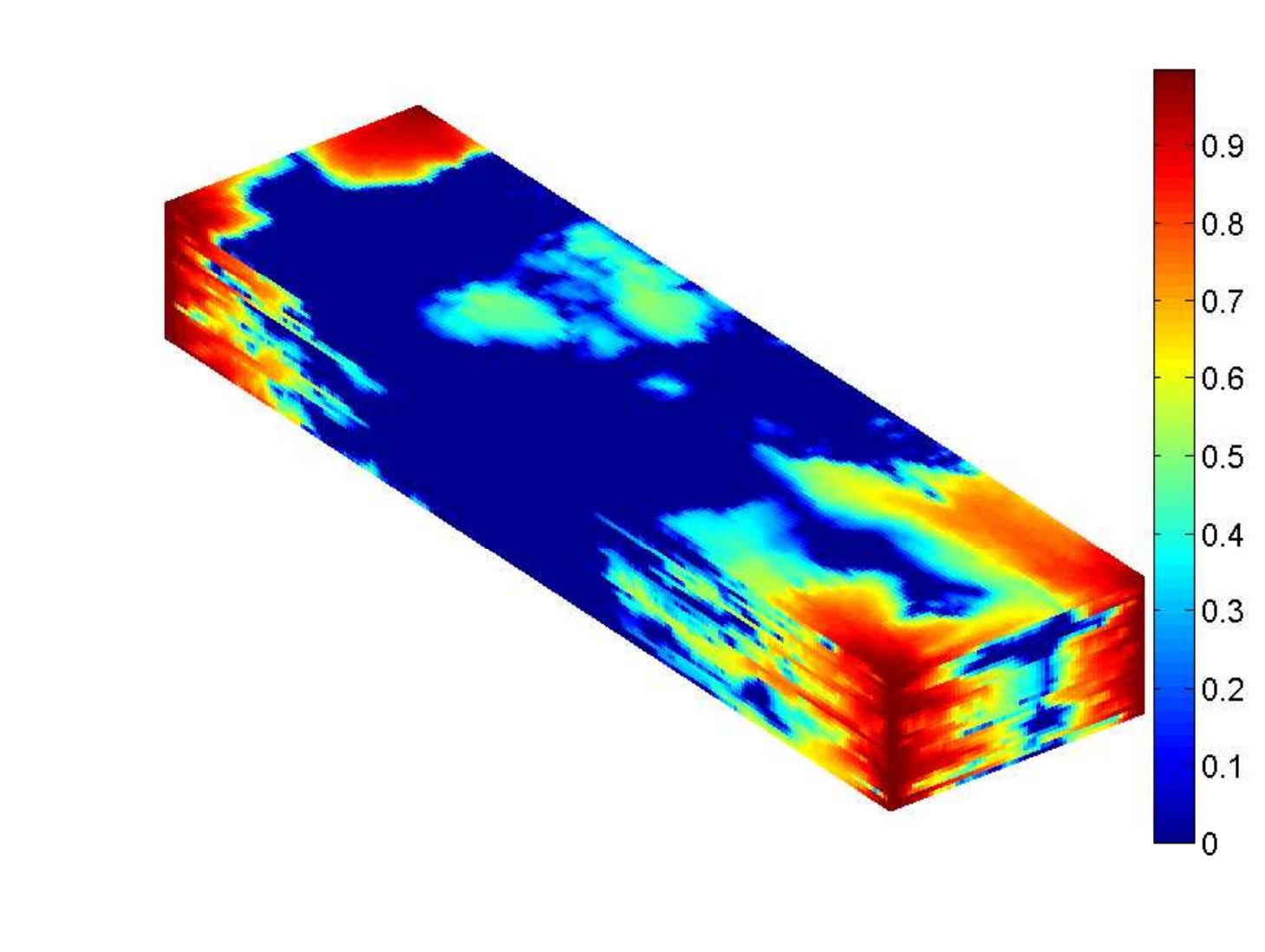}}
	\subfigure[Saturation profile with MMMFEM at $t = 1000$, $n=10$, $Nb=2$]{
		\includegraphics[trim={0 1cm 0  0.8cm},clip,width=4in]{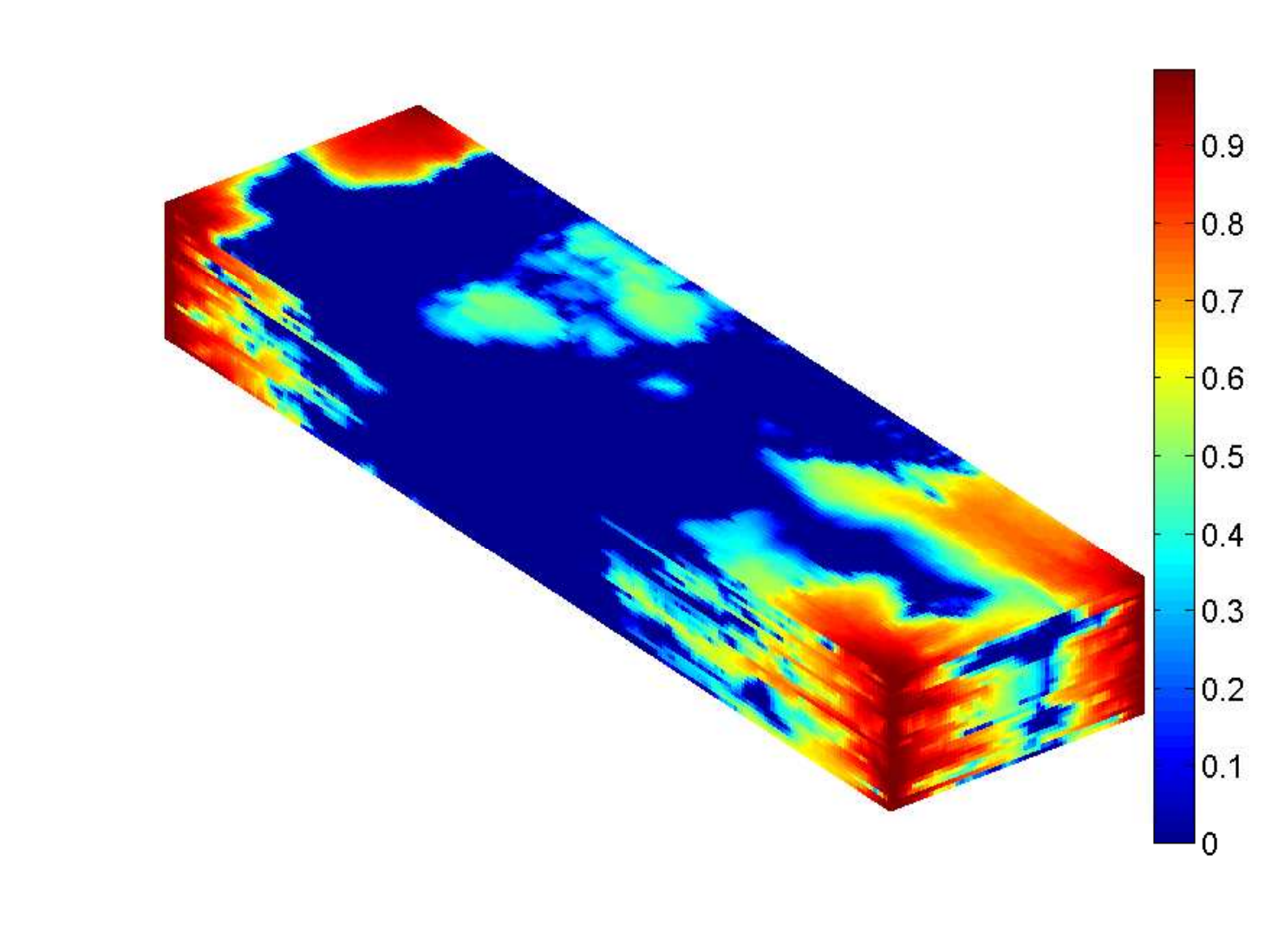}}
	\subfigure[Saturation profile with MMMFEM at $t = 1000$, $n=5$, $Nb=4$]{
		\includegraphics[trim={0 1cm 0  0.8cm},clip,width=4in]{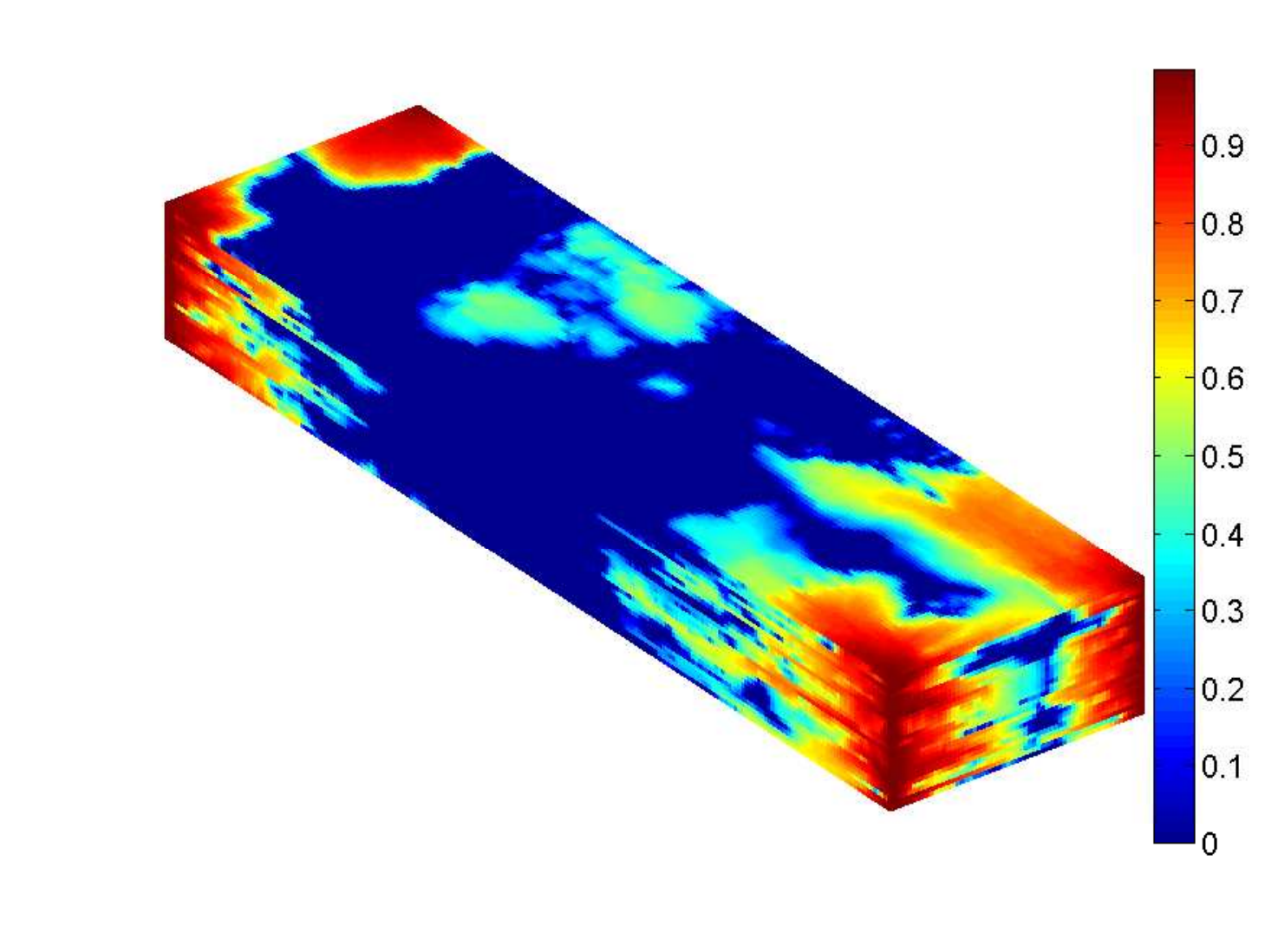}}	
	\caption{Saturation profile at $t = 1000$ for fine-scale solution (top figure), multiscale mortar solution with $n=10$ and $Nb=2$ (middle figure), and multiscale mortar solution with $n=5$  and $Nb=4$ (bottom figure).}
	\label{fig:sat30} 
\end{figure}

The $L^2$ error of the saturation is presented in
Figure \ref{fig:satfirst30}. We can observe that for each case and 
all time instants, the error is less than $10\%$ which clearly demonstrates
the accuracy of the proposed method. Using more basis functions and increase
coarse grid size will definitely reduce the error. We plot the water cut
comparison for model 2 in Figure \ref{fig:pcfirst30}. The lines that 
represent the MMMFEM solution is very close to the reference line.
All lines are almost identical at later times.
Figure \ref{fig:sat30} shows the saturation profiles at $t=1000$ with 
different methods, we observe the results obtained
using the MMMFEM look exactly the same as the fine-scale reference solution.
We report  the accuracy and efficiency comparisons between MMMFEM and fine-scale simulations for model 2 in Table \ref{ta:m2}.
We find that for all MMMFEM solutions, the average $L^2$ error 
of the saturation is less than $5\%$.
The computational time of the MMMFEM decreased by about $60\%$
to $70\%$ compared with the reference solution.

\begin{table}[H]
	\centering \begin{tabular}{|c|c|c|c|c|c|c|}\hline
		Simulator & Basis type [$Nb$]& $n$&Dim&$e_S$ & CPU time (hours)  \tabularnewline\hline
		MMFEM&-&-& 1209600&-& 0.93    \tabularnewline\hline
		MMMFEM&$P_0$+multiscale [2]&10&24192&0.0495 &0.29    \tabularnewline\hline
		MMMFEM&$P_1$+multiscale [4]&10&48384&0.0222 &0.34    \tabularnewline\hline	
		MMMFEM&$P_0$+multiscale [2]&5&96768& 0.0246  &0.32  \tabularnewline\hline
		MMMFEM&$P_1$+multiscale [4]&5&193536& 0.0069   &0.37  \tabularnewline\hline					
	\end{tabular}
	\caption{Accuracy and efficiency comparisons between MMMFEM and fine-scale simulations, model 2.}
\label{ta:m2}	
\end{table}

\begin{figure}[H]
	\centering
	\includegraphics[trim={0 .2cm 0  0cm},clip,width=4in]{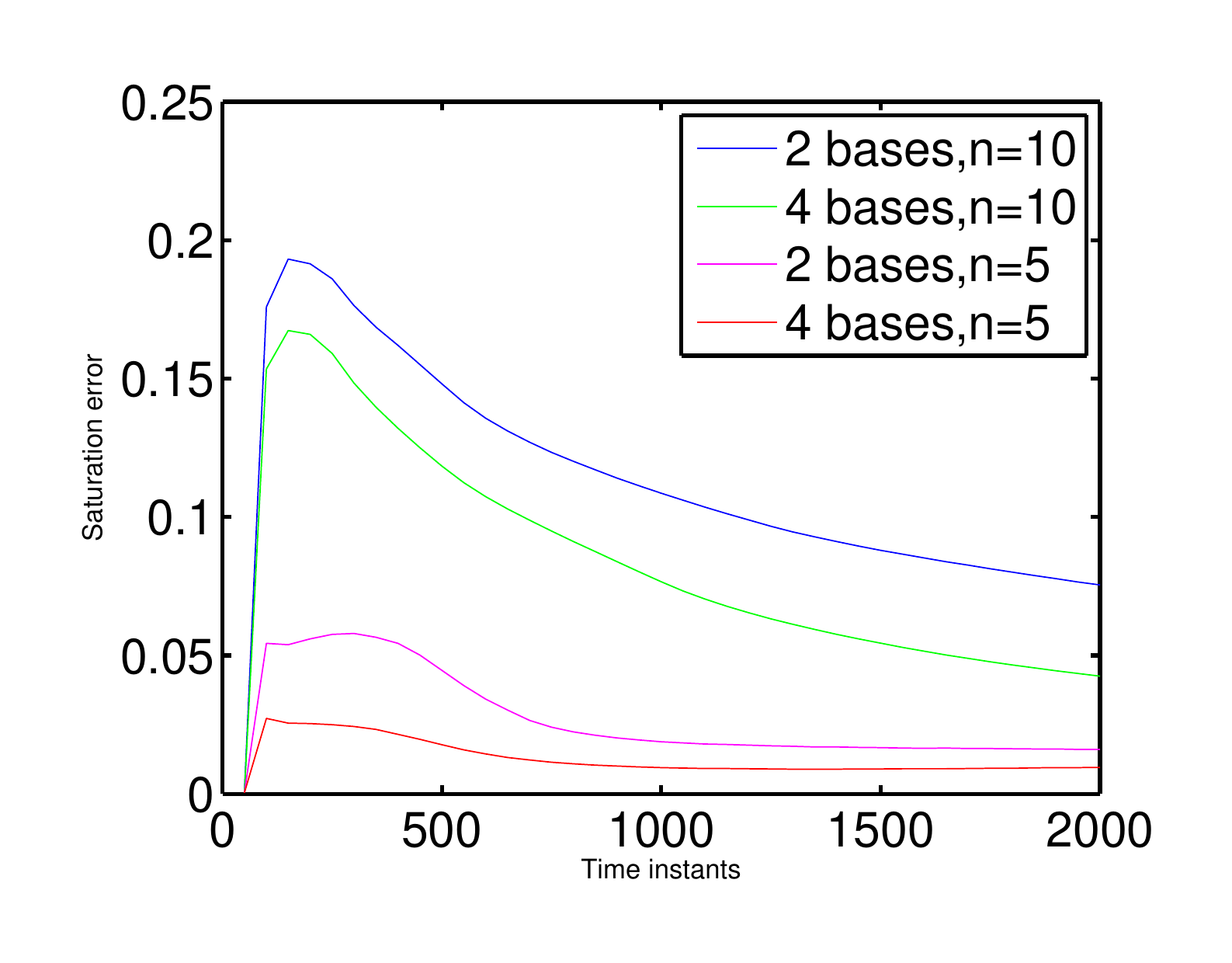}	
	\caption{Saturation error history for model 3.}
	\label{fig:saterrlast50} 
\end{figure}

\begin{figure}[H]
	\centering
	\includegraphics[trim={0 .2cm 0  0cm},clip,width=4in]{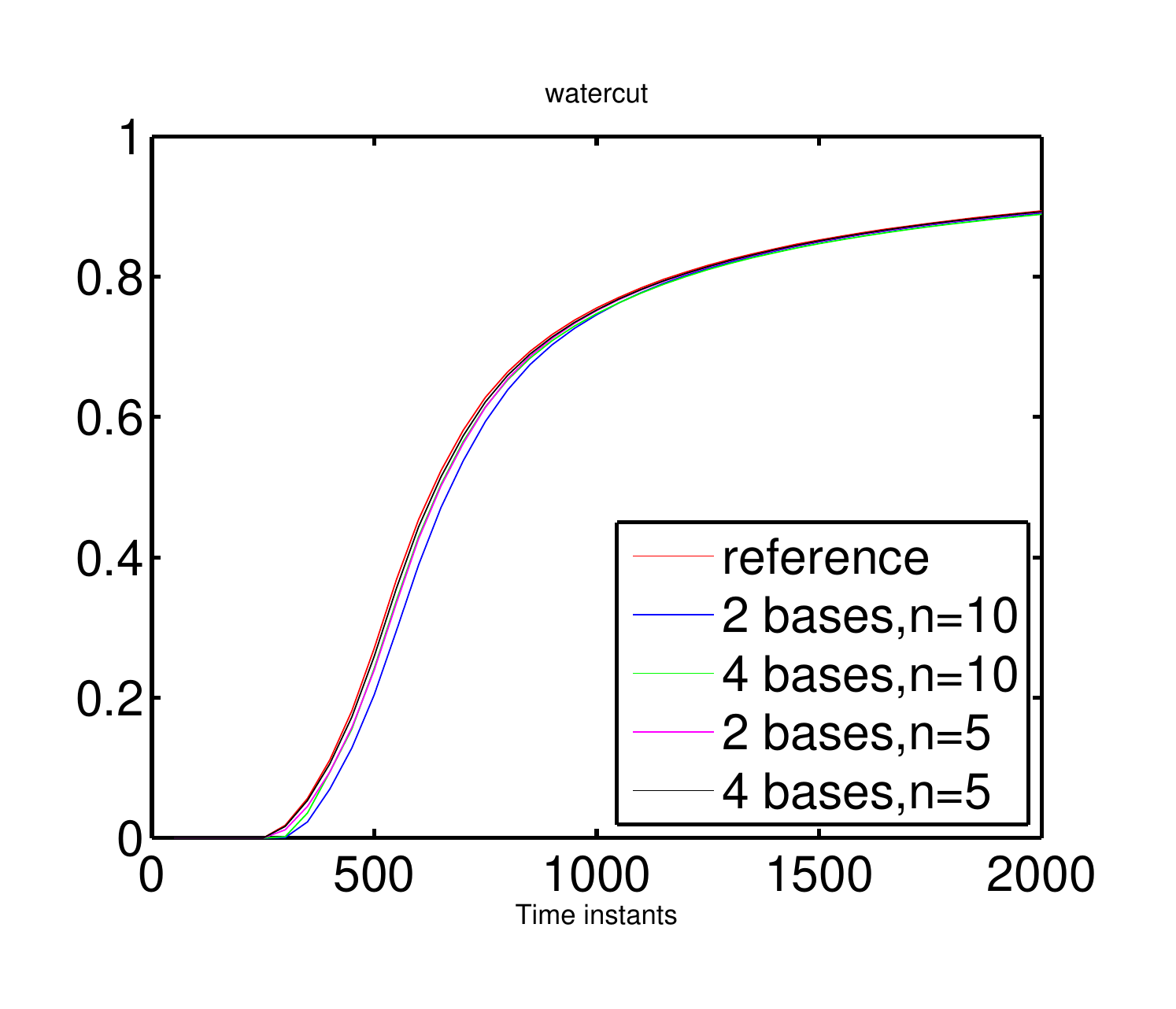}		
	\caption{Watercut comparison for model 3.}
	\label{fig:pclast50} 
\end{figure}

\begin{figure}[H]
	\centering
	\subfigure[Reference saturation profile at $t = 1000$]{
		\includegraphics[trim={0 1cm 0  0.8cm},clip,width=4in]{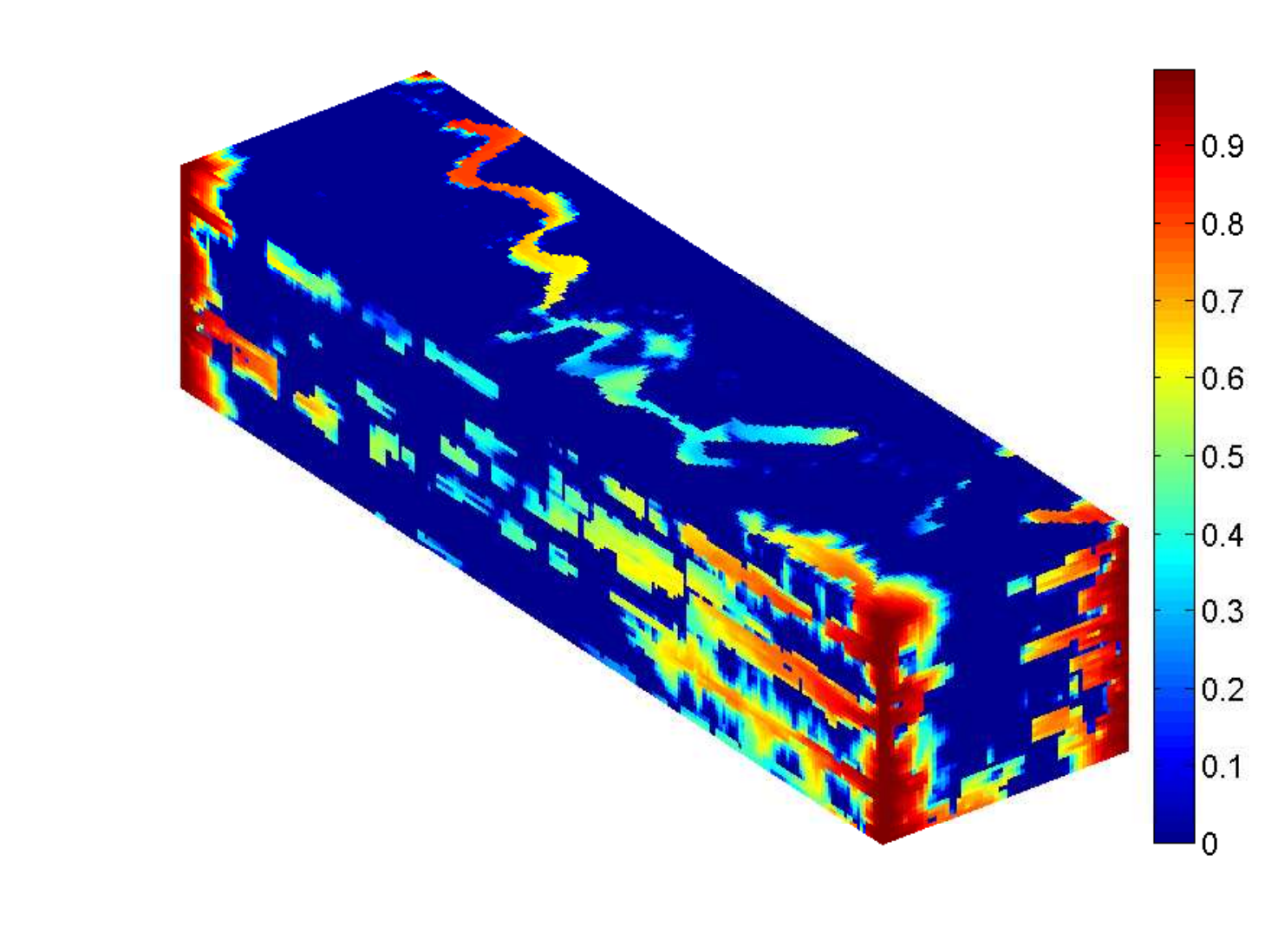}}
	\subfigure[Saturation profile with MMMFEM at $t = 1000$, $n=10$, $Nb=2$]{
		\includegraphics[trim={0 1cm 0  0.8cm},clip,width=4in]{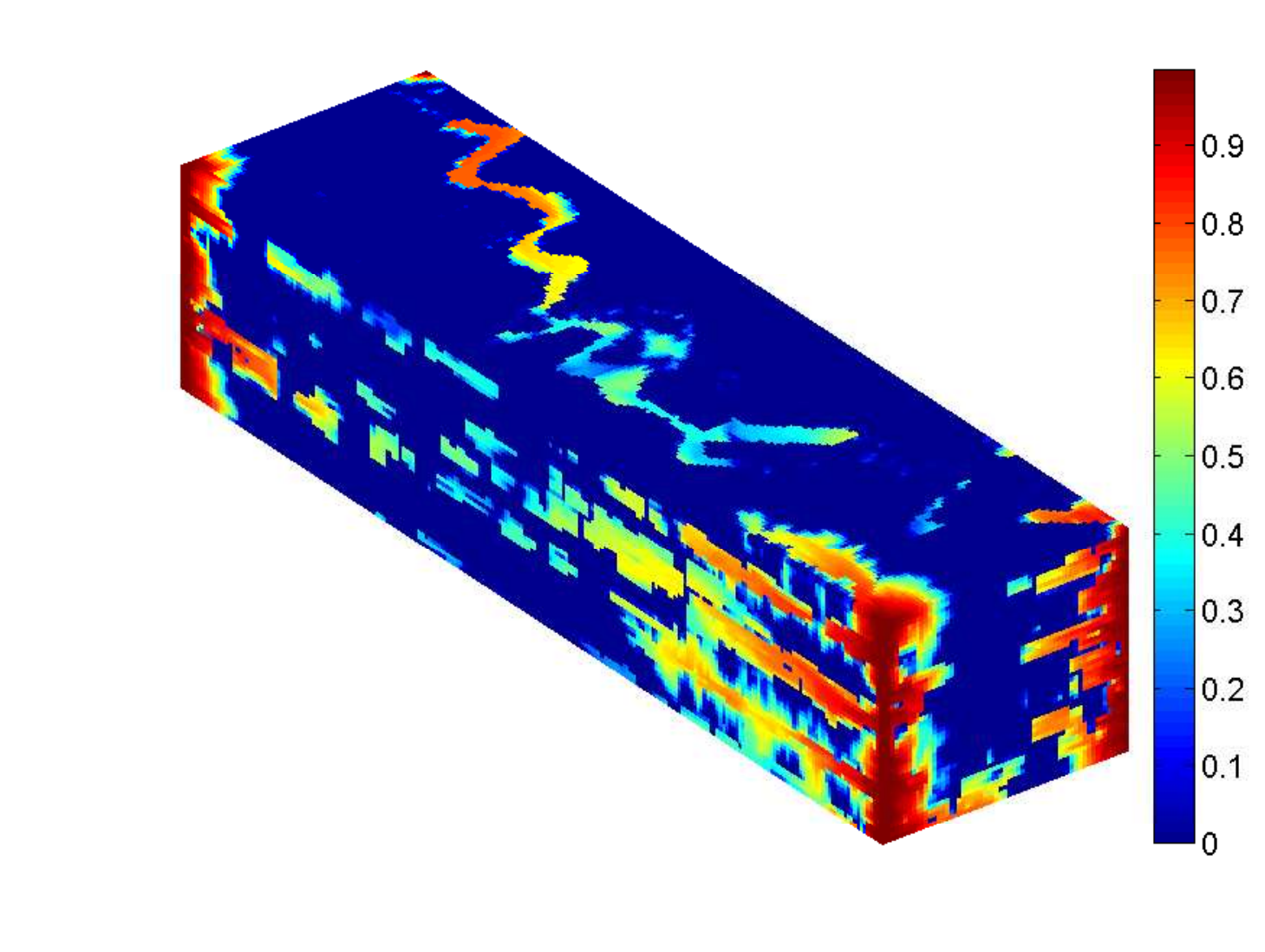}}
	\subfigure[Saturation profile with MMMFEM at $t = 1000$, $n=5$, $Nb=4$]{
		\includegraphics[trim={0 1cm 0  0.8cm},clip,width=4in]{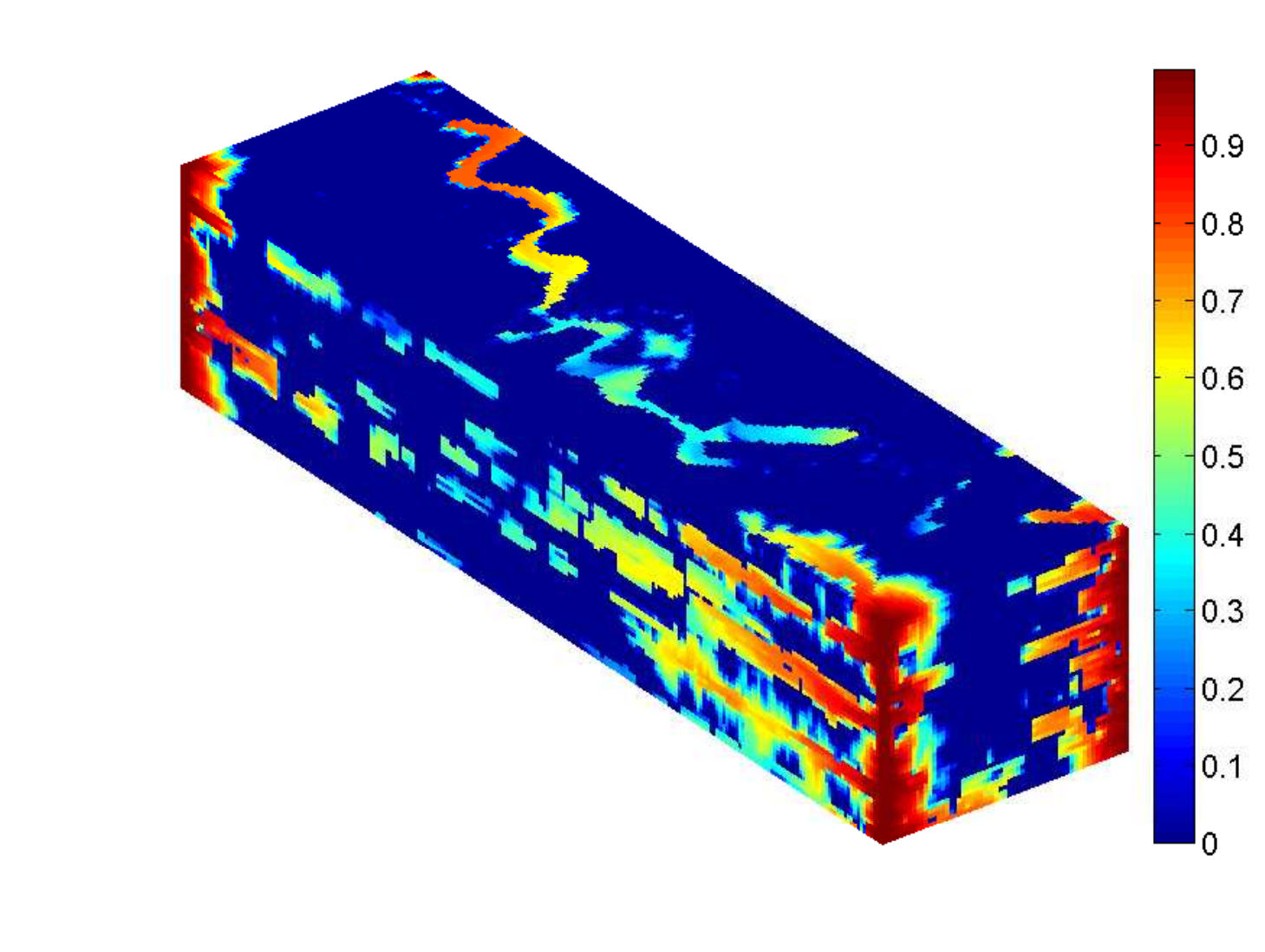}}	
	\caption{Saturation profile at $t = 1000$ for fine-scale solution (top figure), multiscale mortar solution with $n=10$ and $Nb=2$ (middle figure), and multiscale mortar solution with $n=5$ and $Nb=4$ (bottom figure).}
	\label{fig:sat50} 
\end{figure}

For model 3, the $L^2$ error and watercut comparison
are shown in Figures \ref{fig:saterrlast50} and \ref{fig:pclast50} respectively.
We find that if the coarsening ratio is 10, the $L^2$ error is $\mathcal{O}(10^{-1})$, which is not satisfiable. However, if we 
set $n=5$ then the accuracy of the MMMFEM solution will be improved a
lot. The saturation profile comparisons at $t=1000$ 
are provided in Figure \ref{fig:sat50}, we observe that the MMMFEM solution can capture almost all the 
fine-scale features.
We compare the accuracy and efficiency of MMMFEM and fine-scale simulations for model 3, the results are presented in Table \ref{ta:m3}. 
The CPU time savings ranges from $68\%$
to $84\%$, while the average error is about $1\%\sim11\%$. We can expect more computational cost savings for model larger than the test models here. 
We would like to note the relative saturation error at earlier time is larger than later times, this is because the saturation profile has sharp fronts and 
thus the permeability field changes a lot compared with the initial permeability field in some coarse elements.
By comparing the results for model 2 and model 3, 
we find that the performance of the proposed method
is better for smoother media, 
therefore, adaptive 
assign the basis number for different coarse interfaces
based on the media's heterogeneity may save more computational 
cost. We will explore this in our future research.

\begin{table}[H]

	\centering \begin{tabular}{|c|c|c|c|c|c|c|}\hline
		Simulator & Basis type [$Nb$]& $n$&Dim&$e_S$ & CPU time (hours)  \tabularnewline\hline
		MMFEM&-&-& 2007200&-& 3.51    \tabularnewline\hline
		MMMFEM&$P_0$+multiscale [2]&10&40144& 0.1139 & 0.55    \tabularnewline\hline
		MMMFEM&$P_1$+multiscale [4]&10&80288&   0.0837&0.63 \tabularnewline\hline	
		MMMFEM&$P_0$+multiscale [2]&5&160576&0.0270   &0.70  \tabularnewline\hline
		MMMFEM&$P_1$+multiscale [4]&5&321152& 0.0127   &1.13  \tabularnewline\hline									
	\end{tabular}
	\caption{Accuracy and efficiency comparisons between MMMFEM and fine-scale simulations, model 3.}
	\label{ta:m3}	
\end{table}

\begin{figure}[H]
	\centering
	\includegraphics[trim={0 .2cm 0  0cm},clip,width=4in]{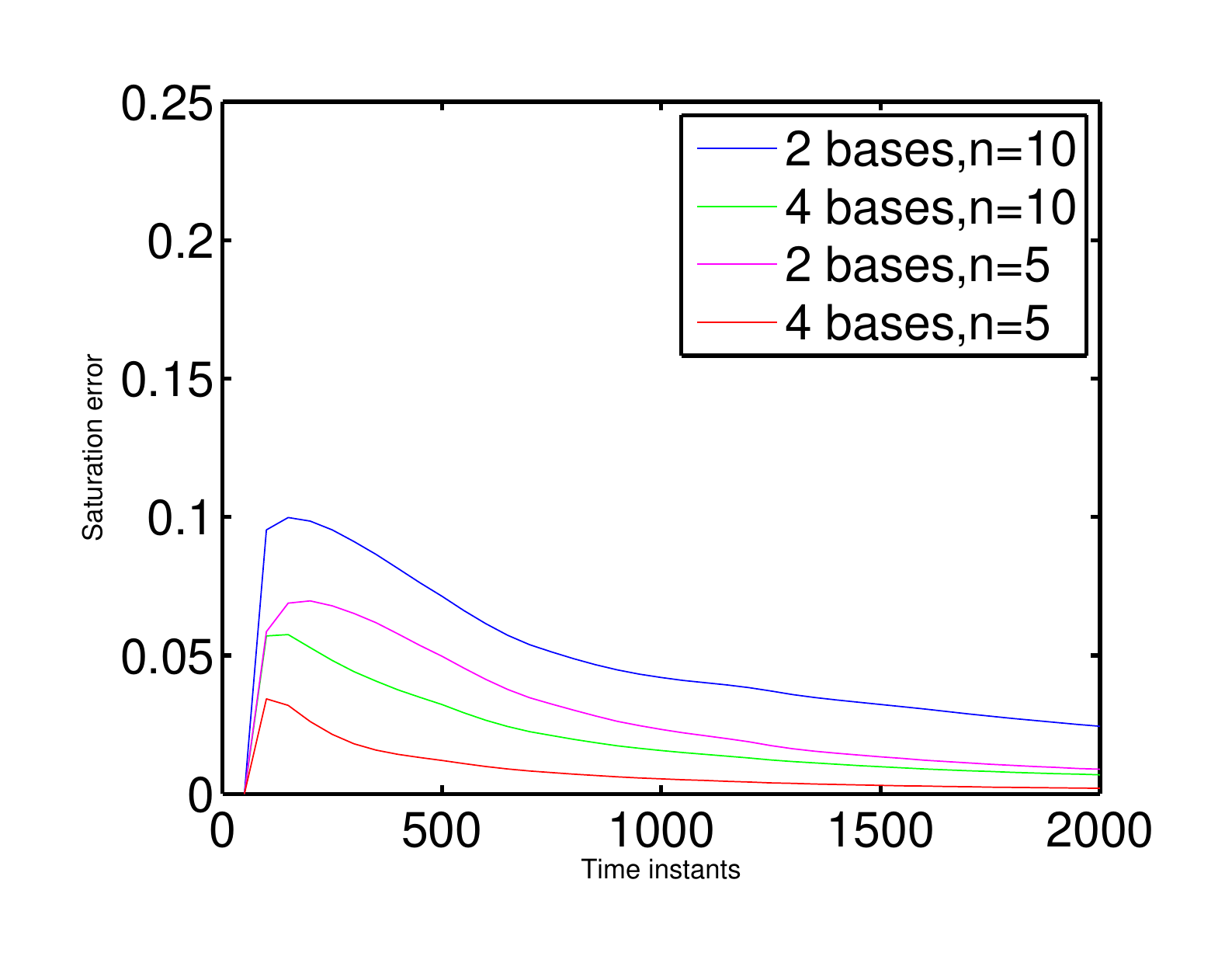}	
	\caption{Saturation error history for model 2, source type II.}
	\label{fig:satfirst30a} 
\end{figure}

\begin{figure}[H]
	\centering
	\includegraphics[trim={0 .2cm 0  0cm},clip,width=4in]{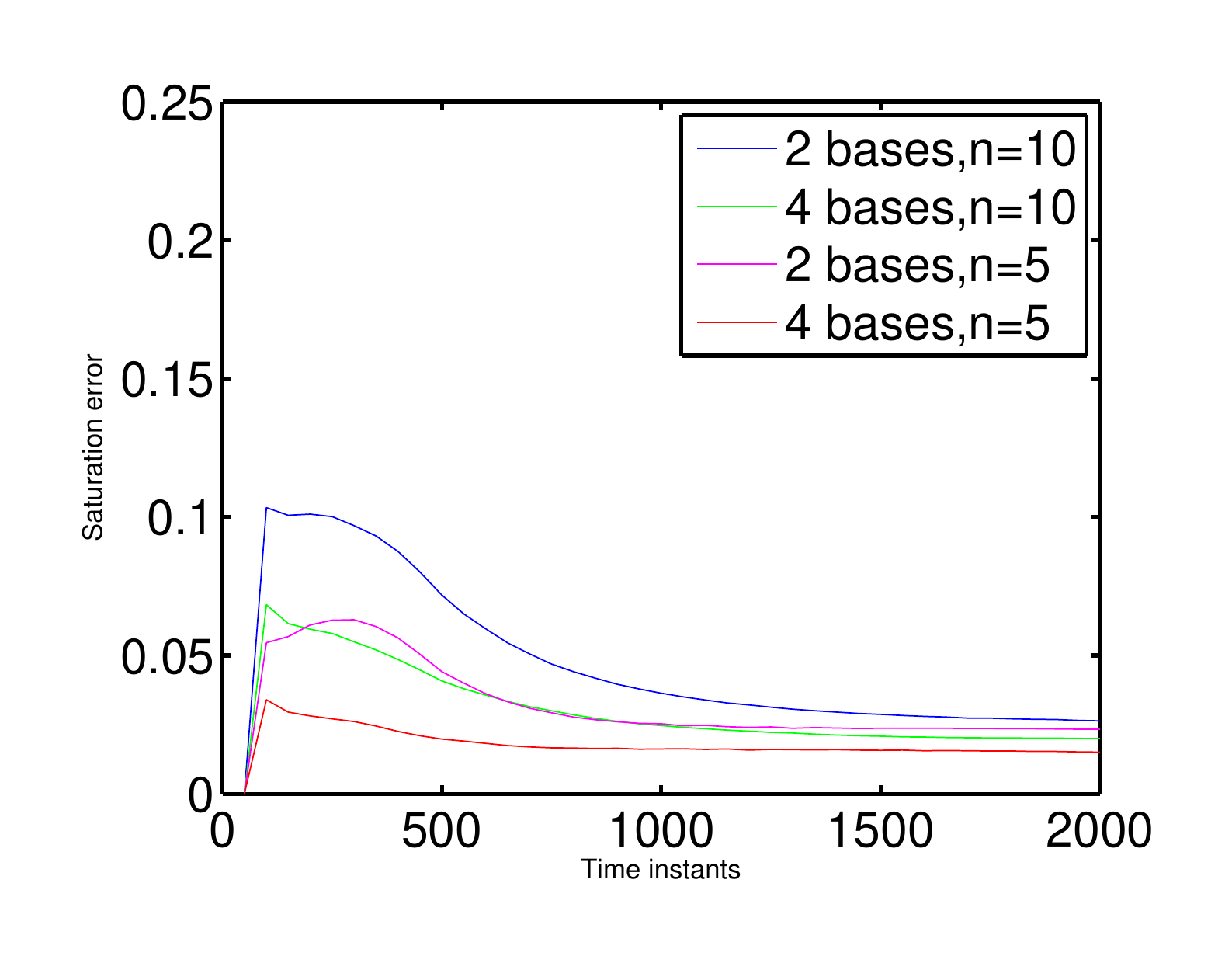}	
	\caption{Saturation error history for model 3, source type II.}
	\label{fig:saterrlast50a} 
\end{figure}

We also considered another type of source (source type II) for the 3D problems, that is, we inject the water in the center of the 
reservoir and place 4 producers in the corners of the reservoir.
The saturation error history for model 2 and 3 are shown in Figures \ref{fig:satfirst30a} and \ref{fig:saterrlast50a} respectively.
Again, we find that the MMMFEM is quite accurate especially for model 3.
Tables \ref{ta:m2a} and \ref{ta:m3a} show the performance comparisons of the MMMFEM and the reference solution.
The average error are all less than $5\%$, similar CPU savings  as previous test cases are observed.

\begin{table}[H]
	\centering \begin{tabular}{|c|c|c|c|c|c|c|}\hline
		Simulator & Basis type [$Nb$]& $n$&Dim&$e_S$ & CPU time (hours)  \tabularnewline\hline
		MMFEM&-&-& 1209600&-&0.93     \tabularnewline\hline
		MMMFEM&$P_0$+multiscale [2]&10&24192&0.0489&0.30  \tabularnewline\hline
		MMMFEM&$P_1$+multiscale [4]&10&48384&0.0204&0.35    \tabularnewline\hline	
		MMMFEM&$P_0$+multiscale [2]&5&96768&0.0289&0.33  \tabularnewline\hline
		MMMFEM&$P_1$+multiscale [4]&5&193536&0.0082&0.38  \tabularnewline\hline					
	\end{tabular}
	\caption{Accuracy and efficiency comparisons between MMMFEM and fine-scale simulations, model 2, source type II.}
\label{ta:m2a}	
\end{table}

\begin{table}[H]
	\centering \begin{tabular}{|c|c|c|c|c|c|c|}\hline
		Simulator & Basis type [$Nb$]& $n$&Dim&$e_S$ & CPU time (hours)  \tabularnewline\hline
		MMFEM&-&-& 2007200&-&3.58    \tabularnewline\hline
		MMMFEM&$P_0$+multiscale [2]&10&40144& 0.0474& 0.57    \tabularnewline\hline
		MMMFEM&$P_1$+multiscale [4]&10&80288&   0.0303&0.64 \tabularnewline\hline	
		MMMFEM&$P_0$+multiscale [2]&5&160576&0.0321&0.69  \tabularnewline\hline
		MMMFEM&$P_1$+multiscale [4]&5&321152&0.0179 &1.09  \tabularnewline\hline									
	\end{tabular}
	\caption{Accuracy and efficiency comparisons between MMMFEM and fine-scale simulations, model 3, source type II.}
\label{ta:m3a}	
\end{table}

\section{Conclusions}
We have developed a local-global multiscale mortar mixed finite element method for the two-phase flow problems in heterogeneous porous media. We use the multiscale mortar mixed finite element method to solve
the pressure equation, and explicit finite volume method to solve the saturation equation.
We propose a dynamic mortar space that consists of polynomials and multiscale
basis, which is the restriction
of global pressure field obtained at previous time step on the coarse interface. 
Fine scale solution of the pressure is obtained to initialize the simulation.
We provide numerical experiments on some benchmark 2D and 3D heterogeneous models 
 to demonstrate the performance of the proposed method. 
\section*{Acknowledgment}

Shubin Fu wants to thank Bingqi Yi (Sun Yat-sen University) for  providing computing resources.

	\bibliographystyle{plain}
\bibliography{references}

\end{document}